\documentclass[11pt]{amsart}

\usepackage{color,graphicx,amssymb,latexsym,amsfonts,txfonts}

\usepackage{pdfsync}

\usepackage{hyperref}
\hypersetup{
    colorlinks=true,       
    linkcolor=blue,          
    citecolor=blue,        
    filecolor=blue,      
    urlcolor=blue           
}

\begin{document}

\newtheorem{theo}{Theorem}[section]
\newtheorem{prop}[theo]{Proposition}
\newtheorem{coro}[theo]{Corollary}
\newtheorem{lemm}[theo]{Lemma}
\newtheorem{rema}[theo]{Remark}
\newtheorem{example}[theo]{Example}

\newtheorem{claim}{Claim}
\newtheorem{conjecture}{Conjecture}
\theoremstyle{remark}

\newcommand{\Ch}{\widehat{\mathbb C}}

\title[Dihedral extended Schottky groups]{Structural description of dihedral extended Schottky groups and application in study of symmetries of handlebodies
}

\author{Grzegorz Gromadzki}
\address{Faculty of Mathematics, Physics and Informatics, Gda\'nsk University, Poland}
\email{grzegorz.gromadzki@ug.edu.pl}

\author{Ruben A. Hidalgo}
\address{Departamento de Matem\'atica y Estad\'{\i}stica, Universidad de La Frontera, Temuco, Chile}
\email{ruben.hidalgo@ufrontera.cl}

\thanks{R. A. Hidalgo was supported by projects Fondecyt 1190001 and 1220261.
G. Gromadzki was supported by Polish NCN 2015/17/B/ST1/03235 and by Chilean Fondecyt 1190001}

\dedicatory{
Dedicated to J.M.R. Sanjurjo for his 70th anniversary
}


\begin{abstract}
Given a  symmetry $\tau$ of a closed Riemann surface $S$, there exists an extended Kleinian group $K$, whose orientation-preserving half is a Schottky group $\Gamma$ uniformizing $S$, such that $K/\Gamma$ induces $\langle \tau \rangle$; the group $K$  is called an extended Schottky group. A geometrical  structural description, in terms of the Klein-Maskit combination theorems, of both Schottky and extended Schottky groups is well known.
A dihedral extended Schottky group is a group generated by the elements of  two different extended Schottky groups, both with the same orientation-preserving half. Such configuration of groups corresponds to closed Riemann surfaces together with two different symmetries and  the aim of this paper is to provide a geometrical structure of them. This result can be used in study of three dimensional manifolds and as an illustration we give the sharp upper bounds for the total number of  connected components of the locus of fixed points of two and three different symmetries of a handlebody with a Schottky structure.
\end{abstract}

\maketitle

\section{Introduction}
Let $S$ be a given closed Riemann surface of genus $g$.
A symmetry of  $S$ is an anticonformal involution $\tau:S \to S$. Harnack's theorem \cite{Harnack} asserts that
each connected component of its set of fixed points (if non-empty) is a simple loop (called an oval or mirror) and that the total number of such ovals is at most $(g+1)$. In algebraic terms, the symmetry $\tau$ corresponds to a real structure on $S$, and this means that $S$ can be described by a real algebraic curve.

\smallskip
By the uniformization theorem, there  exist finitely generated Kleinian groups $\Gamma$, with a $\Gamma$-invariant connected component ${\Omega}$
of its region of discontinuity whose points have trivial $\Gamma$-stabilizer, such that there is a biholomorphism between $S$ and ${\Omega}/\Gamma$.
These say that the pair $(\Gamma,{\Omega})$ is an uniformization of $S$.  A geometrical structural picture of $\Gamma$, in terms of the Klein-Maskit combination theorem \cite{Maskit:Comb, Maskit:Comb4}, was provided in \cite{Maskit:function2}.

\smallskip
Let us consider a finite collection of symmetries $\tau_{1},\ldots,\tau_{n}$ of $S$. We say that these are reflected by the uniformization $(\Gamma,{\Omega})$ of $S$ if there exists an extended Kleinian group $K$ containing $\Gamma$ as a normal subgroup and such that $K/\Gamma$ induces the group $H=\langle \tau_{1},\ldots,\tau_{n}\rangle$. In this setting, one would like to provide a geometrical description of $K$ (which reflects the action of  $H$). In the case when $\Gamma$ is a Fuchsian group and ${\Omega}$ is the upper half-plane (a Fuchsian uniformization), this description is well known (and it comes from the structure of NEC groups).  In this paper, we extend these  to the case of Schottky uniformizations.

\smallskip
A Schottky group of rank $g$ is a purely loxodromic Kleinian group $\Gamma$, with non-empty region of discontinuity $\Omega \subset \widehat{\mathbb C}$, isomorphic to the free group of rank $g$. In this case, $\Omega$ is connected and $\Omega/\Gamma$ is a closed Riemann surface of genus $g$; we say that $(\Gamma,\Omega)$ is a Schottky uniformization of $\Omega/\Gamma$.  Koebe's retrosection theorem \cite{Bers,Koebe} asserts that every closed Riemann surface admits a Schottky uniformization. By the planarity theorem \cite{Maskit:planarity}, Schottky groups provide the lowest uniformizations of closed Riemann surfaces. A geometrical picture of Schottky groups is well known (we recall it in Section \ref{Sec:Schottky}).

\smallskip
An extended Schottky group of rank $g$ is an extended Kleinian group $K$ whose orientation-preserving half $\Gamma$ is a Schottky group of rank $g$.  In this case, $\Omega$ is connected (as both $K$ and $\Gamma$ have the same region of discontinuity) and $\Omega/\Gamma$ is a closed Riemann surface admitting a symmetry which is induced by $K$.
 Conversely, if $\tau$ is a symmetry of $S$, then there is an extended Schottky group $K$, with orientation-preserving half a Schottky group $\Gamma$, there is a biholomorphism $\varphi:S \to \Omega/\Gamma$ (i.e., $(\Gamma,\Omega)$
is a Schottky uniformization of $S$) such that $\varphi \langle \tau \rangle \varphi^{-1}=K/\Gamma$ (i.e., $K$ is an uniformization of the pair $(S,\tau)$). In this way, symmetries of $S$ are in correspondence with extended Schottky groups (whose orientation-preserving halves {\ provide Schottky uniformizations of $S$).
 A geometrical structural description of extended Schottky groups, in terms of the Klein-Maskit combination theorems, was provided in  \cite{H-G:ExtendedSchottky}  (we recall this description in Section \ref{Sec:ExtendedSchottky}). In particular, this description asserts that an extended Schottky group $K$
  is constructed using some $\alpha$ reflections, some $\beta$ imaginary reflections, some $\gamma$ pseudo-hyperbolic transformations, some $\delta$ loxodromic transformations and some $\varepsilon$ real Schottky groups (groups generated by a reflection and a Schottky group keeping invariant the circle of fixed points of the reflection); we set $m(K):=\alpha+\beta+\varepsilon$ (see Section \ref{Sec:ExtendedSchottky} for a geometrical interpretation of this value).

\smallskip
Now, let us assume we are given different symmetries $\tau_{1},\ldots, \tau_{n}$, $n \geq 2$, of the surface $S$ and set $H=\langle \tau_{1},\ldots,\tau_{n}\rangle$. In general, it is not true that there is an extended Kleinian group $K$ containing a Schottky group $\Gamma$ (that unifomizes $S$) as a normal subgroup with $K/\Gamma$ inducing $H$. An example of this situation (for $n=2$) was provided in \cite{H-M:imaginary}. The existence of such a group $K$,  is equivalent to the existence of a collection of pairwise disjoint simple loops of $S$ which is invariant under the action of $H$ and which cut-off $S$ into planar subsurfaces (see Section \ref{Sec:lifting}). If $n=2$, then such an extended Kleinian group $K$ is called a dihedral extended Schottky group of rank $g$; it is
generated by two different extended Schottky groups $K_{1}$ and $K_{2}$, both of rank $g$ and with the same orientation-preserving half.
A geometrical structural description of dihedral extended Schottky groups, containing no reflections, was also obtained in \cite{H-M:imaginary}.   In Sections \ref{Sec:basicos} and \ref{ejemploBEDG}, we describe the basic extended dihedral Schottky groups and our first result (Theorem \ref{const0}) states that every dihedral extended Schottky group (which might or not contain reflections) is constructed from the basic ones.

\smallskip

Let us consider a finite collection of extended Schottky groups $K_{1},\ldots, K_{n}$, all of them of the same rank and with the same orientation-preserving half $\Gamma$, and set $m_{j}=m(K_{j})$.
Then one may ask for an upper bound for the value $m_{1}+\cdots +m_{n}$ and to ask for its sharpness. Our next result  (Theorem \ref{sumak=2}) provides an answer for the cases $n=2,3$.
If we consider the handlebody $M_{\Gamma}=({\mathbb H}^{3}\cup \Omega)/\Gamma$ (we say that $\Gamma$ provides a Schottky structure on it), then each $K_{j}$ induces an orientation-reversing self-homeomorphism of $M_{\Gamma}$ of order two (also called a symmetry) with exactly $m_{j}$ connected components of fixed points (note that the symmetry of $M_{\Gamma}$ restricts to an isometry of the hyperbolic manifold ${\mathbb H}^{3}/\Gamma$ and also restricts to a symmetry of the Riemann surface $\Omega/\Gamma$). The geometrical picture of $K_{j}$ asserts that $m_{j} \leq g+1$ (similar as was for the case of surfaces). In this setting,  the above result provides sharp upper bounds for the total number of components of fixed points of two and three symmetries of a handlebody with a Schottky structure (Corollary \ref{suma}).
This, in particular, asserts that there is no handlebody (with a given Schottky structure) with two different maximal symmetries (those with the maximal number of components) which is not the case for Riemann surfaces (hyperelliptic ones may have two different maximal symmetries).

\subsection*{Notations}
In general, we will use symbols such as $K$ and $\widehat{K}$ to denote either a (extended) Kleinian group or a dihedral extended Schottky group.  For Schottky groups we will use the symbol $\Gamma$.
If $K$ is an extended Kleinian group, then $K^{+}$ will denote its orientation-preserving half, being  its index two subgroup consisting of its M\"obius transformations. The symbol $\Omega$ will denote the region of discontinuity of a Kleinian group $K$ and by $M_{K}=({\mathbb H}^{3} \cup \Omega)/K$ the associated $3$-orbifold.

\section{Preliminaries}\label{Sec:1}
In this section, we recall some facts we will need. More details on these topics may be found, for instance, in \cite{Maskit:book,MT}.

\subsection{Extended Kleinian groups}
M\"obius transformations are classified into {\it parabolic}, {\it loxodromic} (including hyperbolic) and {\it elliptic}. Similarly, extended M\"obius transformations, being the composition of a M\"obius transformations with the  complex conjugation,
are classified into {\it pseudo-parabolic} (the square is parabolic), {\it pseudo-hyperbolic} (the square is
hyperbolic), {\it pseudo-elliptic} (the square is elliptic), {\it reflections} (of order two admitting a circle
of fixed points on $\widehat{\mathbb C}$) and {\it imaginary reflections} (of order two and having no fixed
points on $\widehat{\mathbb C}$). Each imaginary reflection has exactly one fixed point
in ${\mathbb H}^{3}$ and this point determines such imaginary reflection uniquely.
We denote by $\widehat{\mathbb M}$ the group consisting of the M\"obius and extended M\"obius transformations and by ${\mathbb M} \cong {\rm PSL}_{2}({\mathbb C})$ its index two subgroup of M\"obius transformations.
If $K$ is a subgroup of $\widehat{\mathbb M}$, then $K^{+}=K \cap {\mathbb M}$ is its {\it canonical orientation-preserving subgroup}.

\smallskip
A  {\it Kleinian group} (respectively, an {\it extended Kleinian group}) is a discrete subgroup $K$ of ${\mathbb M}$ (respectively,  a discrete subgroup of $\widehat{\mathbb M}$ with $K \neq K^{+}$).
Its {\it region of discontinuity} is the open subset $\Omega$ of $\widehat{\mathbb C}$ (it might be empty) composed by the points on which it acts discontinuously.
Note that $K$ and $K^{+}$  have the same region of discontinuity.

\smallskip
 A {\it function group} is a pair $(K,{\Omega})$, where $K$ is a finitely generated Kleinian group and ${\Omega}$ is a $K$-invariant connected component of its region of discontinuity. In the case that every point in ${\Omega}$ has trivial $K$-stabilizer. The quotient  ${\Omega}/K$ is an analytically finite Riemann surface (a closed Riemann surface with some  finite set of points removed).

\subsection{Klein-Maskit's combination theorems}
In the general theory of abstract groups, amalgamated free products and HNN-extensions are the basic tools that permits to decompose a group into a simplest ones. At the level of (extended) Kleinian groups such a decomposition is known as the Klein-Maskit combination theorem.

\begin{theo}[Klein-Maskit's combination theorem \cite{Maskit:Comb, Maskit:Comb4}]\label{KMC}
\mbox{}
$(1)$ (Free products) Let $K_{j}$ be a (extended) Kleinian group with region of discontinuity $\Omega_{j}$, for $j=1,2$. Let ${\mathcal F}_{j}$ be a fundamental domain for $K_{j}$ and assume that there is simple closed loop $\Sigma$, contained in the interior of  ${\mathcal F}_{1} \cap {\mathcal F}_{2}$, bounding two discs $D_{1}$ and $D_{2}$, so that, for $j \in \{1,2\}$, $\Sigma \cup D_{j} \subset  \Omega_{3-j}$ is precisely invariant under the identity in $K_{3-j}$. Then
\begin{itemize}
\item[({i})] $K = \langle K_1, K_2\rangle$ is a (extended) Kleinian group with fundamental domain ${\mathcal F}_{1} \cap {\mathcal F}_{2}$ and $K$ is the free product of $K_{1}$ and $K_{2}$
    \item[({ii})] every finite order element in $K$ is conjugated in $K$ to a finite order element of either $K_{1}$ or $K_{2}$ and
        \item[({iii})] if both $K_{1}$ and $K_{2}$ are geometrically finite, then $K$ is so.
\end{itemize}
\smallskip
\noindent
$(2)$ (HNN-extensions) Let $K_{0}$ be a (extended) Kleinian group with region of discontinuity $\Omega$, and let ${\mathcal F}_{0}$ be a fundamental domain for $K_{0}$. Assume that there are two pairwise disjoint simple closed loops $\Sigma_{1}$ and $\Sigma_{2}$, both of them contained in the interior of  ${\mathcal F}_{0}$, so that $\Sigma_{j}$ bounds a disc $D_{j}$ such that $(\Sigma_{1} \cup D_{1}) \cap (\Sigma_{2} \cup D_{2})=\emptyset$ and that $\Sigma_{j} \cup D_{j} \subset  \Omega$ is precisely invariant under the identity in $K_{0}$. Let $T$ be either a loxodromic or a pseudo-hyperbolic transformation so that $T(\Sigma_{1})=\Sigma_{2}$ and $T(D_{1}) \cap D_{2}=\emptyset$. Then

\begin{itemize}
\item[({i})]  $K = \langle K_{0}, T \rangle$ is a (extended) Kleinian group with fundamental domain ${\mathcal F}_{1} \cap (D_{1} \cup D_{2})^{c}$ and $K$ is the HNN-extension of $K_{0}$ by the cyclic group $\langle T \rangle$, \item[({ii})] every finite order element of $K$ is conjugated in $K$ to a finite order element of $K_{0}$ and
\item[({iii})] if $K_{0}$ is geometrically finite, then $K$ is so. \hfill $\square$
\end{itemize}

\end{theo}

\subsection{Kleinian $3$-manifolds and their automorphisms}
The group $\widehat{\mathbb M}$ can also be viewed, by the Poincar\'e
extension theorem, as the group of hyperbolic isometries of the hyperbolic space ${\mathbb H}^{3}$ (in this case,
${\mathbb M}$ is the group of orientation-preserving ones).
If $K$ is a Kleinian group with region of discontinuity $\Omega$, then

\medskip
({i}) $M_{K}=({\mathbb H}^{3} \cup \Omega)/K$ is  a
$3$-dimensional orientable orbifold,

({ii}) its interior $M^{0}_{K}={\mathbb H}^{3}/K$ has a hyperbolic structure and

({iii}) its conformal boundary $S_{K}=\Omega/K$ has a natural conformal
structure. \\[2mm] If, moreover, $K$ is torsion free, then $M_{K}$ and $M^{0}_{K}$ are orientable $3$-manifolds and $S_{K}$ is a
Riemann surface; we say that $M_{K}$ is a {\it Kleinian $3$-manifold} and that $M_{K}$ and $S_{K}$ are {\it uniformized by $K$}. Now,
if $K$ is an extended Kleinian group, then the $3$-orbifold $M_{K^{+}}$ admits
the  orientation-reversing homeomorphism $\hat{\tau}:M_{K^{+}} \to M_{K^{+}}$ of order two induced by $K$ and
$M_{K^{+}}/\langle \hat{\tau} \rangle=({\mathbb H}^{3} \cup \Omega)/K$.

\smallskip
Let $\Gamma$ be a torsion free Kleinian group, so $M_{\Gamma}=({\mathbb H}^{3} \cup \Omega)/\Gamma$ is a Kleinian $3$-manifold.
 An {\it automorphism} of $M_{\Gamma}$ is a self-homeomorphism whose restriction to its interior $M_{\Gamma}^{0}$ is a hyperbolic
isometry. An orientation-preserving automorphism is called a {\it conformal automorphism} and an
orientation-reversing one an {\it anticonformal automorphism}. A {\it symmetry} of $M_{\Gamma}$ is an anticonformal
involution. We denote by ${\rm Aut }(M_{\Gamma})$ the group of automorphisms of $M_{\Gamma}$ and by ${\rm Aut}^{+}(M_{\Gamma})$ the subgroup
of conformal automorphisms. If $\pi^{0}:{\mathbb H}^{3} \to M_{\Gamma}^{0}$ is a universal covering induced by $\Gamma$, then it
extends to a universal covering $\pi:{\mathbb H}^{3} \cup \Omega \to M_{\Gamma}$ with $\Gamma$ as the group of Deck
transformations. If $H \subset {\rm Aut }(M_{\Gamma})$ is a finite group and we lift it to the universal covering  space
${\mathbb H}^{3}$ under $\pi^{0}$, then we obtain an (extended) Kleinian group $K$ containing $\Gamma$
as a normal subgroup of finite index such that $H=K/\Gamma$. The group $H$ contains orientation-reversing automorphisms  if and only if
$K$ is an extended Kleinian group.

\subsection{Schottky groups (geometrical picture)}\label{Sec:Schottky}   The {\it Schottky group of rank $0$} is just the trivial group. A {\it Schottky
group of rank $g \geq 1$} is a Kleinian group $\Gamma$ generated by loxodromic
transformations $A_1,\ldots,A_g$, so that there are $2g$ pairwise disjoint simple loops,
$C_1,C'_1,\ldots,C_g,C'_g$, with a $2g$-connected outside   $D\subset \widehat{\mathbb C}$, where
$A_i(C_i)=C'_i$, and $A_i(D)\cap D=\emptyset$, for $i=1,\ldots,g$.
 Its region of discontinuity is connected and dense, $S_{\Gamma}$ is a closed Riemann surface of genus $g$, $M_{\Gamma}$ is topologically a handlebody of genus $g$ and $M_{\Gamma}^{0}$  carries a geometrically finite complete hyperbolic Riemannian metric with injectivity radius bounded away from zero. In this case, we say that $\Gamma$ provides a Schottky structure on the handlebody.
If $g \geq 2$,  then ${\rm Aut }(M_{\Gamma})$ has order at most $24(g-1)$ and ${\rm Aut }^{+}(M_{\Gamma})$ has
order at most $12(g-1) $ \cite{Z1,Z2}. Each conformal (respectively anticonformal) automorphism of $M_{\Gamma}$
induces a conformal (respectively anticonformal) automorphism of the conformal boundary $S_{\Gamma}$ and the
later determines the former due to the Poincar\'e extension theorem.
As a consequence of Koebe's retrosection theorem \cite{Bers, Koebe}, every closed
Riemann surface is isomorphic to $S_{\Gamma}$ for a suitable Schottky group $\Gamma$.

\begin{rema} \rm
(1) A Schottky group of rank $g$ can be equivalently defined as ({i}) a purely loxodromic Kleinian group of the second kind isomorphic to a free of rank $g$ \cite{Maskit:Schottky groups} or ({ii}) as a purely loxodromic geometrically finite Kleinian group isomorphic to a free of rank $g$ (since a free group cannot be the fundamental group of a closed hyperbolic  $3$-manifold).\\[2mm]
(2) The interior of a topological handlebody carries many different types of complete hyperbolic structures. Those, which are geometrically finite and with injectivity radius bounded away from zero,  are exactly the ones provided by Schottky groups.\\[2mm]
(3) Let $M$ be a topological handlebody of genus $g$ and $H$ be a finite group of homeomorphisms of $M$. As $M$ is a compression body,
there are: ({i}) a (extended) Kleinian group $K$, containing as a finite index normal subgroup a Schottky group $\Gamma$ of rank $g$, and ({ii})
an orientation-preserving homeomorphism $f:M \to M_{\Gamma}$, with $f H f^{-1} = K/\Gamma$ \cite{Z2}.
\hfill $\square$
\end{rema}

\subsection{A geometrical structure description of extended Schottky groups}\label{Sec:ExtendedSchottky}
An {\it extended Schottky group of rank $g$} is an extended Kleinian group whose orientation-preserving half
is a Schottky group of rank $g$.  In the case that it does not contains reflections ({\it Klein-Schottky groups}) a  geometrical structure description was provided in \cite{H-M:imaginary}.  A geometric  structural description of all extended Schottky groups, in terms of the Klein-Maskit combination theorems, is as follows.

\begin{theo}[\cite{H-G:ExtendedSchottky}]\label{maintheo}
\mbox{}
An extended Schottky group is the free product $($in
the Klein-Maskit combination theorem sense$)$ of the following kind of groups:
\begin{itemize}
\item[({i})] cyclic groups generated by reflections,
\item[({ii})] cyclic groups
generated by imaginary reflections,
\item[({iii})] cyclic groups generated by
pseudo-hyperbolic,
\item[({iv})] cyclic groups generated by loxodromic
transformations, and
\item[({v})] real Schottky groups $($that is groups generated
by a reflection and a Schottky group keeping invariant the corresponding circle of
its fixed points$)$.
\end{itemize}

\noindent
Conversely,  a group constructed (by the Klein-Maskit theorem) using $\alpha$ groups of type
$({i})$, $\beta$ groups of type $({ii})$, $\gamma$ groups of type ({iii}), $\delta$ groups of type ({iv})  and $\varepsilon$ groups of
type ({v}), is an extended Schottky group if and only if
$\alpha + \beta + \gamma + \varepsilon>0$.
If, in addition, the $\varepsilon$ real Schottky groups above have
the ranks $r_{1},\ldots, r_{\varepsilon} \geq 1$, then it has rank
$g=\alpha + \beta +2(\gamma + \delta) + \varepsilon -1  +r_1 + \ldots + r_\varepsilon.$   \hfill $\square$
\end{theo}

\smallskip
We will say that an extended Schottky group $K$, constructed using $\alpha$ groups of type
$({i})$, $\beta$ groups of type $({ii})$, $\gamma$ groups of type ({iii}), $\delta$ groups of type ({iv})  and $\varepsilon$ groups of
type ({v}) signature $(\alpha,\beta,\gamma,\delta,\varepsilon, \{r_{1},\ldots,r_{\varepsilon}\})$. If $\varepsilon=0$, then we use the notation
$(\alpha,\beta,\gamma,\delta,0, \emptyset)$. We also set $m(K):=\alpha+\beta+\epsilon$, which is the number of connected components of the locus of fixed points of the induced symmetry on the handlebody $M_{\Gamma}$, where $\Gamma$ is the index two orientation-preserving half of $K$.

\begin{coro}\label{fijos de simetria}
If $K$ is an extended Schottky group of signature $(\alpha,\beta,\gamma,\delta,\varepsilon, \{r_{1},\ldots,r_{\varepsilon}\})$ and
$\Gamma$ is its orientation-preserving half, then $K$ induces a  symmetry on the handlebody $M_{\Gamma}$ whose connected components of fixed points consist of
$\alpha$  two dimensional closed discs, $\beta $ isolated points,  and $\varepsilon$ two dimensional non-simply connected compact surfaces.
  \hfill $\square$  \end{coro}

\begin{rema} \rm
Let $K$ be an extended Schottky group of rank $g$ and let $\Gamma$ be its orientation-preserving half. The group $K$ induces a symmetry $\hat{\tau}$ on the handlebody $M_{\Gamma}$. The above facts assert the following observations: \\[2mm]
(1) $\hat{\tau}$ has no fixed points if and only if $K$ has signature $(0,0,\gamma,\delta,0, \emptyset)$, where $\delta>0$
 and $g=2(\gamma+\delta)-1$. In particular, on an even genus handlebody every symmetry must have fixed points.\\[2mm]
(2)  If $n_{0}$ is the number of isolated fixed points, $n_{1}$ the number of total ovals in the conformal boundary
and $n_{2}$ the number of two-dimensional connected components of the set of fixed points of $\hat{\tau}$, then
$n_{0}+n_{1}, n_{0}+n_{2} \in \{0,1,\ldots, g+1\}.$
\hfill $\square$ \end{rema}

\begin{rema}[A connection to real structures] \rm
The space of ${\rm PSL}_{2}({\mathbb C})$-equivalence classes of Schottky representations of $F_{g}$ is the marked Schottky space ${\mathcal MS}_{g}$, which is  known to be isomorphic to the quasiconformal deformation space of a Schottky group of rank $g$ (see, for instance, \cite{Bers,B,Nag}). It is  a non-compact, connected and non-simply connected complex manifold of dimension $3(g-1)$ (a domain of holomorphy of ${\mathbb C}^{3g-3}$ \cite{Nag}) whose (discrete) group of holomorphic automorphisms is isomorphic to ${\rm Out}(F_{g})$ \cite{Earle}.  A real structure on ${\mathcal MS}_{g}$ is an order two anti-holomorphic automorphism $\tau:{\mathcal MS}_{g} \to {\mathcal MS}_{g}$. A fixed point of $\tau$ is called a real point and the set of fixed points  ${\rm Fix}(\tau)$ is its real locus. If the real locus is non-empty, then each of its connected components is a real manifold of dimension $3(g-1)$. Some generalities on real structures on complex manifolds can be found, for instance, in \cite[Section 2]{CC}.  Harnack's problem \cite{Krasnov}  concerns with the number and description of the connected components of the real locus of a real structure. Also, a natural question is about the number of real structures, up to conjugation, that ${\mathcal MS}_{g}$ has.
These questions where answered in \cite{HS}. For instance, in ${\mathcal MS}_{2}$ there are exactly $4$ non-conjugated real structures: one without real points, one with exactly eight real components, one with three and one with four. If $g \geq 3$, then there are exactly $T_{g}+1$ real structures, up to conjugation, where $T_{g}$ is the number of conjugacy classes of elements of order two in ${\rm Out}(F_{g})$.
In the same paper, it was also observed that a real point of a real structure of ${\mathcal MS}_{g}$ can be identified with an extended Schottky group of rank $g$ and that each irreducible component is a real analytic embedding of the quasiconformal deformation space of some extended Schottky group of rank $g$.
If two different real structures of ${\mathcal MS}_{g}$, with non-empty real loci, have a common real point (in which case, these two real structures are necessarily conjugated \cite{HS}), then such a common real point can be identified with a dihedral    extended Schottky groups of rank $g$. \hfill $\square$
\end{rema}

\subsection{A lifting criteria}\label{Sec:lifting}
Below we recall a simple criterion for lifting automorphisms which we will need in Section \ref{Sec:teo2}. This is a direct consequence of the Equivariant Loop Theorem \cite{M-Y}, whose proof is  based on minimal surfaces, that is, surfaces that minimize locally the area. In \cite{Hidalgo:Auto} there is provided a proof whose arguments is proper to (planar) Kleinian groups.

\begin{theo}\cite{Hidalgo:Auto, M-Y}\label{lifting}
Let $(K,{\Omega})$ be a torsion free function group such that $S={\Omega}/K$ is a closed Riemann surface of genus $g \geq 2$.
Let $P:{\Omega} \to S$ be a regular covering with $K$ as its deck group. Let
$H$ be a group of (conformal and anticonformal) automorphism of $S$. Then $H$ lifts with respect to the above regular planar covering if and only if
there is a collection $\mathcal F$ of pairwise disjoint simple loops on $S$ such that:
\begin{itemize}
\item[({i})] $\mathcal F$ defines the regular covering $P:{\Omega} \to S$; and
\item[({ii})] $\mathcal F$ is invariant under the action of $H$.   \hfill $\square$
\end{itemize}
\end{theo}

\begin{rema}\rm
(1) In the above, that $\mathcal F$ defines the regular covering $P:{\Omega} \to S$ means that the covering map is, up to isomorphisms, the smallest covering with the property that the loops in ${\mathcal F}$ lift to loops.\\[2mm]
(2) If $K$ is a Schottky group (so ${\Omega}$ is its region of discontinuity), then condition ({i}) is to say that ${\mathcal F}$ cut-off ${\Omega}/K$ into a collection of planar surfaces.
\hfill $\square$\end{rema}

\subsection{A counting formula}
Let $K$ be an extended Kleinian group containing a Schottky group $\Gamma$ of rank $g$ as a normal
subgroup of finite index (the last trivially holds if  $g \geq 1$). Let  us
denote by $\theta:K \to H=K/\Gamma$,  the canonical projection. If $\hat{\tau} \in H$ is an involution
which is the $\theta$-image of an extended M\"obius transformation, then ${{K}}=\theta^{-1}(\langle
\hat{\tau} \rangle)$ is an extended Schottky group whose orientation-preserving half is $\Gamma$.
Let $(\alpha,\beta,\gamma,\delta,\varepsilon,\{r_{1},\ldots,r_{\varepsilon}\})$ be the signature of ${{K}}$. Some of the values of the signature
can be computed from $\theta$.  Before, we provide some necessary definitions.

\smallskip
A {\it complete set of symmetries of } $K$ is a maximal collection
${\mathcal C}= \{c_{i} \in K: i \in I\}$ of anticonformal involutions (i.e. reflections and
imaginary reflections) which are non-conjugate in $K$( we shall refer to its elements as to canonical symmetries).  As $K$ is geometrically finite (as it is a finite extension of a Schottky group), ${\mathcal C}$ is finite.

\smallskip
 For each $i \in I$ we set $I(i) \subset I$ defined by
those $j \in I$ so that $\theta(c_i)$ and $\theta(c_{j})$ are conjugate in $H$ (in particular, $i \in I(i)$).
Note that it may happen that for $j \in I(i)$, $c_{j}$ can be imaginary reflection even if $c_i$
is a reflection and  viceversa (this occurs when $\theta(c_i)$ is a symmetry of $M_{\Gamma}$ whose locus of fixed points has isolated
points and also two-dimensional components). We set by $J(i)$ the subset of $ I(i)$
defined by those $j$ for which  $c_j$ is an imaginary reflection.
We also set by $F(i) \subset I(i) \setminus J(i)$ for those $j$ for which $c_{j}$ has finite centralizer in
$K$ and  set $E(i)=I(i) \setminus (J(i)\cup F(i))$. As ${{K}}$ has finite index in
$K$, a reflection $c \in K$ has an infinite centralizer ${\rm C}( K,c)$  in
$K$ if and only if it has an infinite centralizer in ${{K}}$.  The crucial for the present paper is the following generalization of a certain formula proved in \cite{gg} for Fuchsian group, to extended Schottky groups obtained in \cite{H-G:ExtendedSchottky}

\begin{theo}\label{teo1}
Let $K$ be an extended Kleinian group, containing a Schottky group $\Gamma$ as a finite index normal
subgroup, and let ${\mathcal C}= \{c_{i}: i \in I\}$ be a complete set of its symmetries.
Let $\theta:K \to H=K/\Gamma$ be the canonical projection and, for $i \in I$, set
${{K}}_{i}=\theta^{-1}(\langle \theta(c_i) \rangle)=\langle \Gamma, c_i \rangle$. Then ${{K}}_{i}$ is an extended Schottky group with signature
$(\alpha,\beta,\gamma,\delta,\varepsilon,\{r_{1},\ldots,r_{\varepsilon}\})$, where
$$\alpha=\sum_{j \in F(i)}
[ {\rm C}(H, \theta(c_j)): \theta({\rm C}(K, c_j))],\;\;
\beta=\sum_{j \in J(i)} [ {\rm C}(H, \theta(c_j)):
\theta({\rm C}(K, c_j))], \;\;
\varepsilon=\sum_{j \in E(i)}
[ {\rm C}(H, \theta(c_j)): \theta({\rm C}(K, c_j))].$$   \hfill $\square$
\end{theo}

\section{Dihedral extended Schottky groups}
A  {\it dihedral extended Schottky group of rank $g$} is an extended Kleinian group generated by two different extended Schottky groups, both with the same Schottky group of rank $g$ as the orientation-preserving half. In this section we provide a geometrical description of these groups in terms of the Klein-Maskit combination theorems.

\subsection{The DESG-criteria}\label{Sec:criterio}
We first proceed to describe a simple criteria which permits us to decide when an extended Kleinian group is a dihedral extended Schottky group.

\begin{prop}[DESG-criteria]\label{necesario}
An extended Kleinian group $K$ is a dihedral extended Schottky
group if and only  there is a surjective homomorphism
$\varphi  :K \to {\rm D}_{p}=\langle a,b: a^{2}=b^{2}=(ab)^{p}=1\rangle $ (for some positive integer $p>1$) whose
kernel is a Schottky group $\Gamma$ and $\varphi  (K^{+})=\langle ab\rangle$.
\end{prop}
\begin{proof}
Let $K$ be a dihedral extended Schottky group, say generated by two different extended Schottky groups ${{K}}_{1}$ and ${{K}}_{2}$, both with the same Schottky group $\Gamma$ as orientation-preserving half. If $\Omega$ is the region of discontinuity of $K$ (which is the same for ${{K}}_{j}$ and $\Gamma$), then ${{K}}_{j}$ induces a symmetry $\hat{\tau}_{j}$ on the Riemann surface $S=\Omega/\Gamma$. If $p>1$ is the order of $\hat{\tau}_{1}\hat{\tau}_{2}$, then
there exists a surjective homomorphism $\varphi:K \to {\rm D}_{p}=\langle \hat{\tau}_{1},\hat{\tau}_{2}\rangle$ such that ${{K}}_{j}=\varphi^{-1}(\langle \hat{\tau}_{j}\rangle)$, $\ker(\varphi)=\Gamma$ and $\varphi(K^{+})=\langle \hat{\tau}_{1}\hat{\tau}_{2} \rangle$.
In the other direction, assume $K$ is an extended Kleinian group and there is a surjective homomorphism $\varphi  :K \to {\rm D}_{p}=\langle a,b: a^{2}=b^{2}=(ab)^{p}=1\rangle$ as required.
As $\varphi  (K^{+})=\langle ab\rangle$, then
${{K}}_{1}=\varphi  ^{-1}(\langle a \rangle)$ and ${{K}}_{2}=\varphi  ^{-1}(\langle b \rangle)$ are extended Kleinian groups. As  $L_{j}$  contains $\Gamma$ as an index two subgroup, it is an
extended Schottky group (where $\Gamma$ is its index two preserving half). Since $K$ is generated by ${{K}}_{1}$ and ${{K}}_{2}$, it is a dihedral extended Schottky
group.
\end{proof}

\begin{rema} \rm
If $K$ is a dihedral extended Schottky group, then its limit set is totally disconnected (as $K$ contains a Schottky group of finite index) and contains no parabolic transformations. Let $N$ be a finite index subgroup of the dihedral group $K$, which is torsion free. We claim that $N$ is a Schottky group. In fact, this is a
a consequence of the classification of function groups \cite{Maskit:function2,Maskit:function3,Maskit:function4}, which states that every finitely generated Kleinian group, without parabolic nor elliptic elements and with totally disconnected limit set, is a Schottky group.
\hfill $\square$
\end{rema}

\subsection{Some examples of dihedral extended Schottky groups}\label{Sec:basicos}
Let $K$ be an extended Kleinian group constructed (see Figure \ref{Fig1}), by use of the Klein-Maskit combination theorem, as a free product of $\alpha$ cyclic groups generated by a reflection, $\beta$ cyclic groups generated by an imaginary reflections, $\gamma$ cyclic groups generated by a loxodromic transformation, $\delta$ cyclic groups generated by a pseudo-hyperbolic, $\rho$ cyclic groups generated by an elliptic transformation of finite order and $\eta$ groups generated by an elliptic transformation of finite order and a loxodromic transformation commuting with the elliptic one and $\kappa$ groups generated by an elliptic transformation of finite order and a pseudo-hyperbolic transformation both sharing their fixed points (so the glide reflection conjugates the elliptic into its inverse).
If $\Omega$ is the region of discontinuity of $K$, then
\begin{enumerate}
\item[({i})] $\Omega/K$ is the connected sum of $\beta+2\gamma+2\delta+2\kappa+2\eta$ real projective planes (if $\beta+\delta+\kappa>0$) or a genus $\gamma+\eta$ orientable surface (if $\beta=\delta=\kappa=0$),
 with exactly $\alpha$ boundary components and with $2\rho$ conical points in the interior (in pairs of the same cone order); and
\item[({ii})] $\Omega/K^{+}$ is a closed Riemann surface of genus $\widetilde{g}=\alpha+\beta+2(\gamma+\delta+\kappa+\eta)-1$ with $4\rho$ conical points (in quadruples with the same cone order).
\end{enumerate}
As a consequence of the DESG-criteria (Proposition \ref{necesario}),  $K$ will be a dihedral extended Schottky group if and only if either
\begin{enumerate}
\item[({i})] $\alpha+\beta+\delta+\kappa \geq 2$ or
\item[({ii})]  $\alpha+\beta+\delta+\kappa=1$ and $\gamma+\eta \geq 1$.
\end{enumerate}
In that proposition, one may take $p$ as the least common multiple of the orders of the $\rho+\eta+\kappa$ elliptic generators and the epimorphism $\varphi$ sending the $\gamma+\eta$ loxodromic transformations to the identity or powers of $ab$, the reflections, imaginary reflections and pseudo-hyperbolics to either $a$ or $b$ (or their conjugate) and the elliptic elements  to powers of $ab$ (preserving the order). These examples will correspond to the ones constructed in Theorem \ref{const0} (by using only the groups of types ({i})-({vii}) in there stated).

\begin{figure}[h]
\begin{center}

\includegraphics[width=8cm,keepaspectratio=true]{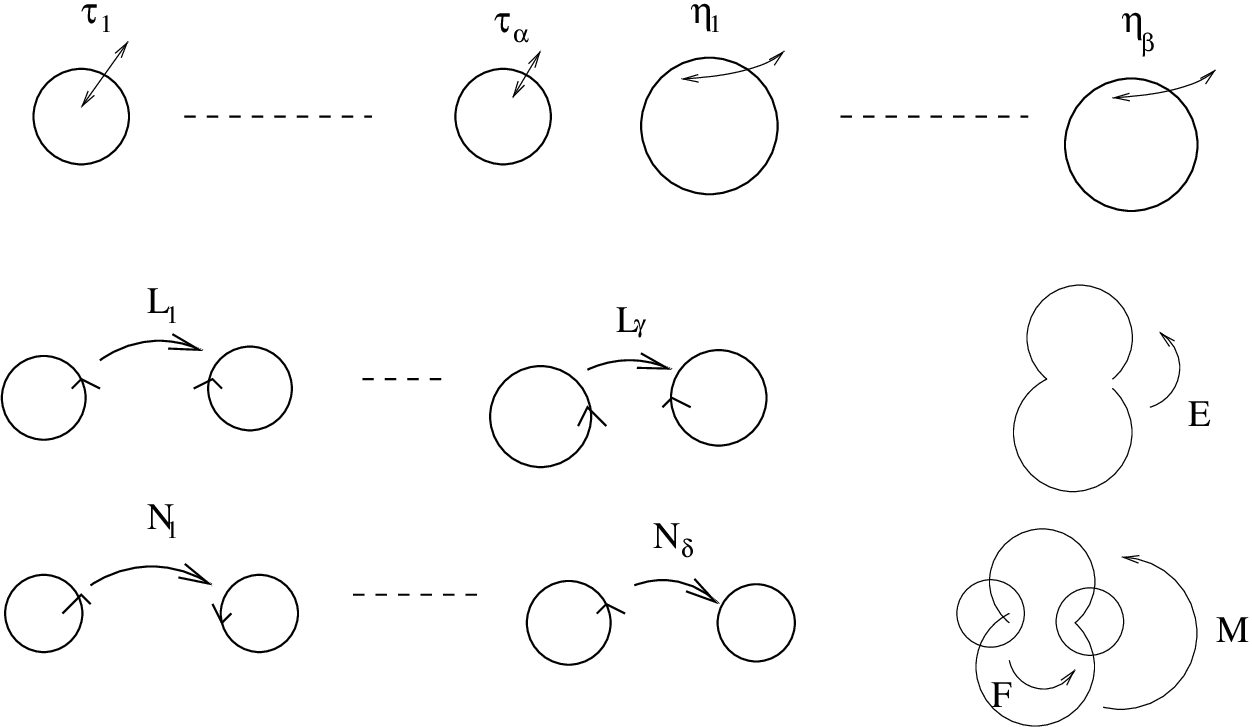}
\caption{Examples of dihedral extended Schottky groups: $\tau_{i}$ are reflections, $\eta_{i}$ are imaginary reflections, $L_{i}$ are loxodromics, $N_{i}$ are pseudo-hyperbolics, $E,F$ are elliptics, $M$ is either loxodromic or a pseudo-hyperbolic.}
\label{Fig1}
\end{center}
\end{figure}

\subsection{Basic extended dihedral groups}\label{ejemploBEDG}
Next, we proceed to define the basic extended dihedral Schootky groups, needed in the statement of Theorem \ref{const0}.
Let $C \subset \widehat{\mathbb C}$ be a circle and let $\sigma$ be the reflection on it. Consider a (finite) non-empty collection of pairwise disjoint closed discs ${\Omega}_{j}$, each one with boundary circle $\Sigma_{j}$ being orthogonal to $C$. We denote by ${\rm Int}({\Omega}_{j})$ and ${\rm Ext}({\Omega}_{j})$ the interior and exterior, respectively, of ${\Omega}_{j}$.

\begin{enumerate}
\item Take $2r$ of these circles (it could be $r=0$), say $\Sigma_{1}, \cdots,\Sigma_{2r}$, and loxodromic transformations $L_{1},\ldots,L_{r}$ such that $L_{j}(\Sigma_{j})=\Sigma_{r+j}$, $L_{j}({\rm Int}({\Omega}_{j})) \cap {\rm Int}({\Omega}_{j+r})=\emptyset$ and each $L_{j}$ commuting with $\sigma$. (The group generated by $\sigma$ and $L_{1},\ldots,L_{r}$ is a real Schottky group).

\item Now, for some each others $\Sigma_{i}$, we consider an elliptic transformation $E_{i}$, with both fixed points on $C$, such that $E_{i}({\rm Ext}({\Omega}_{i})) \subset {\rm Int}({\Omega}_{i})$ and  $\sigma E_{i} \sigma=E_{i}^{-1}$ (each $E_{i}\sigma$ is a reflection).

\item For others $\Sigma_{k}$, we consider an elliptic transformation $F_{k}$, with both fixed points on $C$, such that $F_{k}({\rm Ext}({\Omega}_{k})) \subset {\rm Int}({\Omega}_{k})$ and a loxodromic transformation $M_{k}$, such that $M_{k} F_{k}=F_{k}M_{k}$, $M_{k}\sigma=\sigma M_{k}$ and $\sigma F_{k} \sigma=F_{k}^{-1}$ (both fixed points of $F_{k}$ are also the fixed points of $M_{k}$).

\item For others $\Sigma_{s}$ we take an elliptic transformation of order two $D_{s}$ (whose fixed points are not on $C$) keeping invariant $\Sigma_{s}$ (so permuting both discs bounded by it) and commuting with $\sigma$.

\item Finally, for each of  the rest of the circles $\Sigma_{l}$ we consider the reflection $\tau_{l}$ on it (so it commutes with $\sigma$).
\end{enumerate}

The group $K$ generated by $\sigma$ and all of the above transformations is an extended Kleinian group (by the Klein-Maskit combination theorem).
If $\Omega$ is the region of discontinuity of $K$, then $S=\Omega/K$ is a bordered (orbifold) Klein surface, which is orientable if and only if the loxodromic transformations $L_{k}$ keep invariant each of the two discs bounded by $C$. Each elliptic element $E_{i}$ produces two conical points on the border, both with cone order the order of $E_{i}$. Each $D_{s}$ produces a conical point of order two in the interior. Each reflection in (5) produces two conical points of order two in the boundary. In particular, $\Omega/K^{+}$ provides of a compact orientable orbifold with some even number of conical points admitting a symmetry permuting them.

\medskip

Let $p$ be the least common multiple of the orders of all elliptic transformations $E_{i}$, $F_{k}$ and $D_{s}$. The epimorphism $\varphi:K \to D_{p}=\langle a,b: a^{2}=b^{2}=(ab)^{p}=1\rangle$,  sending
the loxodromic transformations $L_{j}$ and $M_{k}$ to the identity, the reflection $\sigma$ to $a$ and the reflections $\tau_{l}$ to $b$, satisfies the conditions of the DESG-criteria (Proposition \ref{necesario}). In particular,
$K$ is a dihedral extended Schottky group; called a {\it basic extended dihedral group} (see Figure \ref{Fig2}).

\begin{figure}[h]
\begin{center}

\includegraphics[width=8cm,keepaspectratio=true]{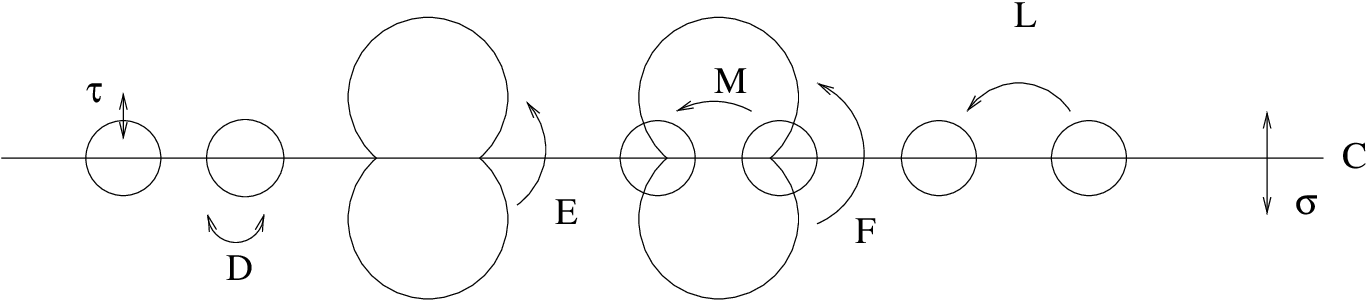}
\caption{Examples of a basic extended dihedral group.}
\label{Fig2}
\end{center}
\end{figure}

\subsection{The structural picture of dihedral extended Schottky groups}\label{Sec:estructura}
In Sections \ref{Sec:basicos} and \ref{ejemploBEDG}, we have constructed (by use of the DESG-criteria) some examples of dihedral extended Schottky groups. Next result asserts that all dihedral extended Schottky groups are obtained, by use of the Klein-Maskit combination theorems, from them.

\begin{theo}\label{const0}
\mbox{}
A dihedral extended Schottky group is the free product (in the sense of the Klein-Maskit combination theorems) of the following groups (see Figures \ref{Fig1} and \ref{Fig2})
\begin{itemize}
\item[({i})] $\alpha$ cyclic groups generated by reflections,
\item[({ii})] $\beta$ cyclic groups generated by imaginary reflections,
\item[({iii})] $\gamma$ cyclic groups generated by loxodromic transformations,
\item[({iv})] $\delta$ cyclic groups generated by pseudo-hyperbolics,
\item[({v})] $\rho$ cyclic groups generated by an elliptic transformation of finite order,
\item[({vi})] $\eta$ groups generated by an elliptic transformation of finite order and a loxodromic transformation commuting with the elliptic one,
\item[({vii})] $\kappa$ groups generated by an elliptic transformation of finite order and a pseudo-hyperbolic transformation both sharing their fixed points (so the glide reflection conjugates the elliptic into its inverse),
\item[({viii})] $\varepsilon$ basic extended dihedral groups $K_{1},\ldots,K_{\varepsilon}$;
\end{itemize}
such that $\alpha+\beta+\delta+\kappa+\varepsilon>0$ and either
$\alpha+\beta +\delta +\kappa+\varepsilon \geq 2$ or
$\alpha+\beta +\delta + \kappa+\varepsilon =1$ and $\gamma+\rho+\eta>0.$   \hfill $\square$
\end{theo}

\begin{rema}\label{observa1} \rm
Let $K$ be a group constructed from the groups ({i})-({vii}) in the above theorem. We define its signature as the tuple $(\alpha,\beta,\gamma,\delta,\rho,\eta,\kappa,\varepsilon)$ and the number $m(K):=\alpha+\beta+\varepsilon$ (which provides the number of orientation-reversing factors used in the construction of $K$). The condition $\alpha+\beta+\delta+\kappa+\varepsilon>0$ is equivalent for $K$  to have orientation-reversing elements, that is, to be an extended Kleinian group. The constructed group $K$ might not be a dihedral extended Schottky group; for it to happen one needs  the conditions ({i}) or ({ii}).

\smallskip
Let $K_{1},\ldots,K_{\varepsilon}$ be the  $\varepsilon$ basic extended dihedral groups used in the construction of $K$, let $\Omega_{i}$ be  the region of discontinuity of $K_{i}$ and let $K^{+}_{i}$ the orientation-preserving half of $K_{i}$.
Each orbifold $\Omega_{i}/K^{+}_{i}$ has genus $g_{i}$, $2r_{i}$ conical points of orders $t_{i1},t_{i1},\ldots, t_{ir_{i}}, t_{ir_{i}}\geq 3$ and $2n_{i}$ conical points of order two. In this case,
the region of discontinuity $\Omega$ of $K$ is connected and the orbifold $\Omega/K^{+}$ has genus
$\widetilde{g}=\alpha+\beta+2(\gamma+\delta+\kappa+\eta)+\varepsilon+g_{1}+\cdots+g_{\varepsilon} -1$,
has $4\rho$ conical points of orders $l_{1},l_{1},l_{1},l_{1},\ldots,l_{\rho},l_{\rho},l_{\rho},l_{\rho}$,
$2(r_{1}+\cdots+r_{\varepsilon})$ conical points, of orders $t_{11}$, $t_{11}$, $t_{12}$,
$t_{12}$, \ldots, $t_{\varepsilon r_{\varepsilon}}$, $t_{\varepsilon r_{\varepsilon}}$ and
$2(n_{1}+\cdots+n_{\varepsilon})$ conical points of order $2$.
If $K$ happens to be a dihedral extended Schottky group (so it satisfies either condition ({i}) or ({ii}) of the previous theorem), then it contains  a Schottky group $\Gamma$ of rank $g$
as a normal subgroup such that $K/\Gamma \cong {\rm D}_{p}$ (where ${\rm D}_{p}$ is the dihedral group of order $2p$ generated by two symmetries of $\Omega/\Gamma$) and (by the Riemann-Hurwitz formula)
$$g=p(\widetilde{g}-1)+1+2p\sum_{i=1}^{\rho}\left(1-\frac{1}{l_{i}}\right)+p\sum_{i=1}^{\varepsilon}\sum_{k=1}^{r_{i}}\left(1-\frac{1}{t_{ik}}\right)+
\frac{p}{2}\sum_{i=1}^{\varepsilon}n_{i}.$$ \hfill $\square$
\end{rema}

\subsection{Counting orientation-reversing factors of two and three extended Schottky groups with same orientation-preserving half}

In the following,  $[x]$ stands for the integer part of the real number $x$.

\begin{theo}\label{sumak=2} Let $g \geq 2$. Then \\[2mm]
$(1)$  Let $\widehat{K}=\langle {{K}}_{1}, {{K}}_{2} \rangle$ be a dihedral extended Schottky group of rank $g$, generated by the extended Schottky groups ${{K}}_{1}$ and ${{K}}_{2}$, with
${{K}}_{1}^{+}={{K}}_{2}^{+}=\Gamma$. Then if
 $m_{j}=m({{K}}_{j})$ and  the dihedral group $H=\widehat{K}/\Gamma$ is of order $2q$, then
 $$m_{1}+m_{2} \leq 2\left[ \frac{g-1}{q}\right]+4.$$
 Moreover, for every integer $q\geq 2$, the above upper bound is sharp for
infinitely many values of $g$.\\[2mm]
$(2)$ Let $\widehat{K}=\langle {{K}}_{1}, {{K}}_{2}, {{K}}_{3} \rangle$ be the extended Kleinian group generated by the extended Schottky groups (of the same rank $g$) ${{K}}_{1}$, ${{K}}_{2}$ and ${{K}}_{3}$, with ${{K}}_{1}^{+}={{K}}_{2}^{+}={{K}}_{3}^{+}=\Gamma$. Then if $H=\widehat{K}/\Gamma$ and, $m_{j}=m(K_{j})$, then
$$m_{1}+m_{2}+m_{3}\leq \left\{ \begin{array}{ll}
5 & \mbox{if} \; g=2.\\
8 & \mbox{if} \; g=3.\\
g+5 &  \mbox{if} \; g \geq 4 \mbox{ and $H \ncong {\mathbb Z}_{2} \times {\rm D}_{r}$ for any $r$.}\\
\dfrac{(r+1)g+5r-1}{r} &  \mbox{if} \; g \geq 4 \mbox{ and $H \cong {\mathbb Z}_{2} \times {\rm D}_{r}$ for some $r$.}\\
\end{array}
\right.
$$
Moreover, the above upper bounds are sharp for $g=2,3$ and for infinitely  many values of $g \geq 4$.   \hfill $\square$
\end{theo}

The above upper bounds are obtained in Section \ref{Sec:pruebasumak=2}. In Section \ref{Sec:4}, we construct explicit examples to see that they are sharp.

\begin{rema} \rm
The possibilities for $H$ in case (2) are provided in Table \ref{tablita}. We may assume
$m_{1} \leq m_{2} \leq m_{3}$.
If $g=2$, then (by (1)) $m_{i}+m_{j} \leq 4$, so the maximal value $m_{1}+m_{2}+m_{3}=5$ can holds for $(m_{1},m_{2},m_{3}) \in \{(1,1,3),(1,2,2)\}$ (the example we construct corresponds to $(1,1,3)$). Similarly, for $g=3$, the maximal value $m_{1}+m_{2}+m_{3}=8$ only holds for  $(m_{1},m_{2},m_{3}) \in \{(2,2,4),(2,3,3)\}$ and this happen for $H \cong {\mathbb Z}_{2}^{3}$ (the example we construct corresponds to $(2,2,4)$).
\hfill $\square$
\end{rema}

\subsection{An application to symmetries on handlebodies}\label{Sec:handlebodies}
Let $\Gamma$ be a Schottky group of rank $g$, with region of discontinuity $\Omega$. In this case, $S_{\Gamma}=\Omega/\Gamma$ is a closed Riemann surface of genus $g$,
$M_{\Gamma}=({\mathbb H}^{3} \cup \Omega)/\Gamma$ is a handlebody of genus $g$ with a Schottky structure provided by $\Gamma$, and
its interior $M_{\Gamma}^{0}={\mathbb H}^{3}/\Gamma$ carries a natural complete hyperbolic structure whose conformal boundary is $S_{\Gamma}$.

Each symmetry $\hat{\tau}$ of $M_{\Gamma}$ induces by restriction a symmetry of $S_{\Gamma}$, which we still denote as $\hat{\tau}$.
By lifting $\langle \hat{\tau} \rangle$ to the universal cover, we obtain an extended Schottky group $G_{\hat{\tau}}$ whose orientation-preserving half is $\Gamma$. As a consequence of the geometrical structure of extended Schottky groups \cite{H-G:ExtendedSchottky}, the locus of fixed points of $\hat{\tau}$ has at most $g+1$  connected components ($\hat{\tau}$ is called maximal if it has $g+1$ connected components of fixed points) and each of such connected components is either ({i}) an isolated point in $M_{\Gamma}^{0}$ or ({ii}) a $2$-dimensional bordered compact surface (which may or not be orientable) whose border is contained in its conformal boundary $S_{\Gamma}$ (see also \cite{Ka-Mc}). The projection of an isolated fixed point of $\hat{\tau}$ produces a point in the orbifold $M_{\Gamma}^{0}/\langle \hat{\tau} \rangle = {\mathbb H}^{3}/{{K}}_{\hat{\tau}}$ admitting a neighborhood which  locally looks like a cone over the projective plane.

\smallskip
Now, assume we are given different symmetries $\hat{\tau}_{1}, \ldots, \hat{\tau}_{n}$ of $M_{\Gamma}$, and set $H=\langle \hat{\tau}_{1}, \ldots, \hat{\tau}_{n} \rangle$. If we lift each of the elements of $H$ to ${\mathbb H}^{3}$, then we obtain an extended Kleinian group $\widehat{K}$ and a canonical surjective homomorphism $\theta:\widehat{K} \to H$ whose kernel is $\Gamma$. The group $\widehat{K}$ is generated by the extended Schottky groups ${{K}}_{\hat{\tau}_{1}}, \ldots, {{K}}_{\hat{\tau}_{n}}$. Each subgroup $\langle {{K}}_{\hat{\tau}_{r}}, {{K}}_{\hat{\tau}_{s}}\rangle$ (where $1 \leq r<s \leq n$) is a dihedral extended Schottky group. Associated to each ${{K}}_{\hat{\tau}_{j}}$ we have associated the integer $m_{j}=m({{K}}_{\hat{\tau}_{j}})$. In general, sharp upper bounds for the sum $m_{1}+\cdots+m_{n}$ are not known. A consequence of Theorem \ref{sumak=2}, we may obtain such sharp upper bounds for the cases $n=2,3$.

\begin{coro}\label{suma}
Let $M$ be a handlebody  of genus $g \geq 2$, with a given Schottky structure. Then \\[2mm]
$(1)$ If $\hat{\tau}_{1}$ and  $\hat{\tau}_{2}$ are two different symmetries of $M$, $q \geq 2$ is the order of $\hat{\tau}_{1}\hat{\tau}_{2}$ and $m_{j}$ is the number of
connected components of fixed points of $\hat{\tau}_{j}$, then
 $$m_{1}+m_{2} \leq 2\left[ \frac{g-1}{q}\right]+4.$$
 Moreover, for every integer $q\geq 2$, the above upper bound is sharp for
infinitely many values of $g$.\\[2mm]
$(2)$ If $\hat{\tau}_{1}, \hat{\tau}_{2}$  and $\hat{\tau}_{3}$ are three different symmetries, $H=\langle \hat{\tau}_{1},\hat{\tau}_{2}, \hat{\tau}_{3}\rangle$ and  $m_{j}$ is
the number of connected components of fixed points of $\hat{\tau}_{j}$,  then
$$m_{1}+m_{2}+m_{3}\leq \left\{ \begin{array}{ll}
5 & \mbox{if} \; g=2.\\
8 & \mbox{if} \; g=3.\\
g+5 &  \mbox{if} \; g \geq 4 \mbox{ and $H \ncong {\mathbb Z}_{2} \times {\rm D}_{r}$ for any $r$.}\\
\dfrac{(r+1)g+5r-1}{r} &  \mbox{if} \; g \geq 4 \mbox{ and $H \cong {\mathbb Z}_{2} \times {\rm D}_{r}$ for some $r$.}\\
\end{array}
\right.
$$
Moreover, the above upper bounds are sharp for $g=2,3$ and for infinite many values of $g \geq 4$.
  \hfill $\square$
\end{coro}

Corollary \ref{suma} asserts the following fact (already observed in \cite{H-M:imaginary} if both symmetries only have isolated fixed points).

\begin{coro}\label{maximales}
Let $\hat{\tau}_{1}$ and $\hat{\tau}_{2}$ be two different symmetries of a handlebody of genus $g \geq 2$, with a Schottky structure, such that $\hat{\tau}_{j}$ has $m_{j}$ connected components of fixed points. Then $m_{1}+m_{2} \leq g+3$. In particular, ({i}) a handlebody of genus $g \geq 2$ admits at most one maximal symmetry and ({ii}) the upper  bound $m_{1}+m_{2}=g+3$ only occurs if $q=2$, that is, when
$\langle \hat{\tau}_{1}, \hat{\tau}_{2}\rangle={\mathbb Z}_{2}^{2}$.
  \hfill $\square$
\end{coro}

\begin{rema}[Connection to symmetries of Riemann surfaces]\label{observa-diferencia} \rm
Let $S$ be a closed Riemann surface.\\
(1) If $\tau_{1}, \tau_{2}$ are two different symmetries of $S$, then in \cite{BCS} it was observed that the total number of ovals of these two symmetries is bounded above by $2(g-1)/q +4$ (if $q$ is odd) and $4g/q+2$ (if $q$ is even), where $q$ is the order of the product of them. If moreover, $q \geq 3$ and $q$ does not divides $g-1$, then the sharp upper bound is $[ 2(g-1)/q]+3$ \cite{Ewa}. We may observe that the upper bounds given in Corollary \ref{suma} are different from the ones obtained for the case of Riemann surfaces. This difference comes from the fact that there are Riemann surfaces admitting two different symmetries which cannot be realized by a dihedral extended Schottky group. The number of ovals of two symmetries is at most $2g+2$, in particular, there must be the possibility for $S$ to admit two different maximal symmetries.
In \cite{Natanzon}, Natanzon proved that if a Riemann surface admits two maximal symmetries, then it is necessarily hyperelliptic. This is again a different situation as that for handlebodies with a Schottky structure (in such a handlebody there is at most one maximal symmetry).\\[1mm]
(2) In \cite{N1} there were obtained a sharp upper bound for the number of ovals for three different symmetries of $S$, this being $2(g+2)$, again different from that of Corollary \ref{suma}.
\hfill $\square$
\end{rema}

\section{Proof of Theorem \ref{const0}}\label{Sec:teo2}
Proposition \ref{necesario}  asserts that conditions ({i}) and ({ii}) of the theorem are necessary for a dihedral extended Schottky group.  We only need to prove that a dihedral extended Schottky group $K$ can be constructed using the basic groups ({i})-({viii}).
Let $K$ be generated by two different extended Schottky groups
${{K}}_{1}$ and ${{K}}_{2}$, such that ${{K}}_{1}^{+}={{K}}_{2}^{+}=\Gamma$ (a Schottky group of rank $g$). Let $K^{+}$ be its index two
orientation-preserving half.  All the groups $\Gamma$, ${{K}}_1$, ${{K}}_2$, $K$ and
$K^{+}$ have the same region of discontinuity $\Omega$. Set $S^+=\Omega/\Gamma$, $S_1=\Omega/{{K}}_1$, $S_2=\Omega/{{K}}_2$, $S_{+}=\Omega/K^{+}$ and
$S=\Omega/K$.

\subsection{Elliptic elements of $K^{+}$}
\begin{prop}  If
$k \in K^{+}$ is an elliptic transformation with one fixed point in $\Omega$, then both fixed points of
$k$ belong to $\Omega$.
\end{prop}
\begin{proof}
As $[K^{+}:\Gamma]<\infty$ and $\Gamma$ contains no  parabolic elements, the same holds for $K^{+}$.
As $K^{+}$ is a geometrically finite function group, the result follows from  \cite{Hidalgo:MEFP}.
\end{proof}

As $K$ is a finite extension of $\Gamma$, and $S^+$ is a closed Riemann surface,
({i}) $S$ is a compact not necessarily orientable orbifold, with a finite number of orbifold points and possible
non-empty boundary, and ({ii}) $S_{+}$ is a closed Riemann surface with some finite number of orbifold points.

\subsection{Orientation-reversing elements of $K$}
Since $\Gamma$ has index 2 in ${{K}}_{i}$, for $i=1, 2$, if $\eta \in {{K}}_{i} \setminus \Gamma$, then
$\eta_{i}^{2} \in \Gamma$ and  $\eta$ is either: ({i}) a reflection or ({ii}) an imaginary reflection or ({iii}) pseudo-hyperbolic.
In particular, ${{K}}_{i}$ induces
a symmetry $\hat{\tau}_{i}$ on $S^+$ with $S_{i}=S^+/\langle\tau_{i}\rangle$ and each orientation-reversing element of $K$ either acts without fixed
points on $\Omega$ or it is a reflection.

\medskip
The group $J=\langle \hat{\tau}_{1}, \hat{\tau}_{2}\rangle$ is isomorphic to the dihedral group ${\rm D}_{p}$, where $p \geq 2$ is the order of $\hat{\tau}_{1}\hat{\tau}_{2}$.
It follows that ({i}) $S_{+}=S^{+}/\langle \hat{\tau}_{1}\hat{\tau}_{2}\rangle$, $S=S^+/J$ and that ({ii}) every orientation-reversing element in $J$ is $J$-conjugate to either $\hat{\tau}_1$ or $\hat{\tau}_2$.
On $S_{+}$ we have an anticonformal involution induced by $J$, preserving the finite set of orbifold points; so the quotient of
$S_{+}$ by it is again $S$.

\smallskip
Let us denote by $\Phi:K\to J$ the canonical surjective homomorphism, with kernel $\Gamma$, which sends the elements of ${{K}}_{i}\setminus \Gamma$ to the symmetry $\hat{\tau}_{i}$. It follows that $\Phi^{-1}(\langle \hat{\tau}_{i}\rangle)={{K}}_{i}$.

\begin{rema}\rm
If $p$ is odd, then $\hat{\tau}_{1}$ and $\hat{\tau}_{2}$ are conjugated in $J$. In particular, ${{K}}_{1}$ and ${{K}}_{2}$
are conjugated in $K$.
\hfill $\square$
\end{rema}

\begin{prop}\label{reversing}
Every orientation-reversing element of $K$ is either pseudo-hyperbolic or an imaginary
reflection or a reflection.
\end{prop}
\begin{proof}
An orientation-reversing element $a \in K$ is a lift of a conjugate of either $\hat{\tau}_1$ or $\hat{\tau}_2$, so $a^2\in \Gamma$. In this way, either ({i}) $a^2=1$, in which case $a$ is either a reflection or an
imaginary reflection, or ({ii}) $a^{2}$ is a loxodromic transformation, in which case it is a pseudo-hyperbolic.
\end{proof}

\subsection{Structure loops and regions}
Let $P:\Omega \to S^{+}$ be a regular planar Schottky covering,  with deck group $\Gamma$. The group $J$ lifts, under $P$, to $K$ and $\Phi(k)P=P k$, for every $k \in K$.
Theorem \ref{lifting} asserts the existence of a $J$-invariant collection of pairwise disjoint loops $\mathcal F$ in $S^{+}$,
dividing $S^{+}$ into genus zero surfaces (since we are dealing with a Schottky covering) and such that each of these loops lifts to a
simple loop on $\Omega$. If $A$ is a connected component of $S^{+} \setminus {\mathcal F}$ and $J_{A}$ is its $J$-stabilizer,
then (as $J$ is a dihedral group) the subgroup $J_{A}$ is either:
\begin{itemize}
\item[({i})] trivial; or
\item[({ii})] a cyclic conformal group generated by a power of $\hat{\tau}_{1}\hat{\tau}_{2}$; or
\item[({iii})] a cyclic group of order two generated by a symmetry (which is conjugated to either $\hat{\tau}_{1}$ or $\hat{\tau}_{2}$); or
\item[({iv})] a dihedral subgroup of $J$.
\end{itemize}

\noindent
In either case ({iii}) or ({iv}),  $A$ is invariant under some symmetry $\hat{\tau} \in J$. If
$\hat{\tau}$ is either ({i}) a reflection containing a loop of fixed points in $A$ or ({ii}) an imaginary reflection, then
we may find a simple loop $\beta \subset A$ which is invariant under $\hat{\tau}$; moreover, if $\hat{\tau}$ is reflection,
then $\beta$ is formed of only fixed points of it. We may add such a loop and its $J$-translated to ${\mathcal
F}$ without destroying the conditions of Theorem \ref{lifting}. In this way, we may also  assume that ${\mathcal F}$ satisfies the following extra property:

\begin{enumerate}
\item[({v})] Every symmetry in $J_{A}$ is necessarily a
reflection whose circle of fixed points is not completely contained in $A$ (it must intersect some boundary loops).
\end{enumerate}

We may also assume ${\mathcal F}$ to be minimal in the sense that there is no proper subcollection of it satisfying ({i})-({v}) above.
The loops in ${\mathcal F}$ are called the {\it base structure loops} and the connected components of $S^{+} \setminus {\mathcal F}$ are called the {\it base structure
regions}. Both collections are $J$-invariants.

\smallskip
Let $\mathcal G$ be the collection of loops on $\Omega$ obtained by the lifting of those in $\mathcal F$; these are called
 the {\it structure loops}.  The connected components
of $\Omega \setminus {\mathcal G}$ will be called the {\it structure regions}.  Both collections, ${\mathcal G}$ and the set of structure regions, are $K$-invariant.

\subsection{Stabilizers of structure regions}
Let $R$ be a structure region and $K_{R}$ be its $K$-stabilizer. As $P:R \to A=P(R)$ is a homeomorphism, the restriction homomorphism $\Phi: K_{R} \to J_{A}$
is an isomorphism. As a consequence, of the above properties on ${\mathcal F}$, we obtain the following fact.

\begin{prop}\label{regstab}
If $R$ is a structure region, then $K_{R}$ is either:
\begin{itemize}
\item[({i})] trivial;
\item[({ii})] a finite cyclic group generated by an elliptic transformation;
\item[({iii})] a cyclic group of order two generated by a reflection
whose circle of fixed points is not completely contained on $R$;
\item[({iv})] a dihedral group generated two reflections, each one
with its circle of fixed points not completely contained on $R$.
\end{itemize}
\end{prop}

\begin{prop} If $R$ is a structure region with $K_{R}=\{I\}$,
then the restriction to $R$ of the projection map from $\Omega$ to $S=S^{+}/J=\Omega/K$ is a homeomorphism
onto its image.
\end{prop}
\begin{proof}
If this is not the case, then there are two different points $x, y \in R$ and some $k \in K$ such that $k(x)=y$. But in this case, $k(R) \cap R \neq \emptyset$. As the collection of structure regions is $K$-invariant, it follows that $k(R)=R$, i.e., $k\in K_{R}=\{I\}$, a contradiction.
\end{proof}

\subsection{Stabilizers of structure loops}
If $\beta \in {\mathcal G}$, then we denote its $K$-stabilizer as $K_{\beta}$.
As $J^{+}=\langle \hat{\tau}_{1}\hat{\tau}_{2}\rangle$ (a cyclic group), the orientation-preserving half of $K_{\beta}$ is either trivial or a finite cyclic
group generated by some elliptic element.
By Proposition \ref{reversing}, an orientation-reversing transformation inside $K$ is
either an imaginary reflection or a reflection or a pseudo-hyperbolic transformation. As the structure loop $\beta$
is contained in $\Omega$, a pseudo-hyperbolic cannot belong to $K_{\beta}$.
In particular, the only orientation-reversing transformations in $K$ that keep invariant
some structure loop can be either an imaginary reflection or a reflection.
Also, $K_{\beta}$ cannot contain two different imaginary reflections (since the product of two distinct imaginary reflections is always
hyperbolic).

\smallskip
The structure loop $\beta$ can be stabilized by a
reflection in two different manners. One is that it fixes it point-wise, that is, $\beta$ is
the circle of fixed points of the reflection. The second one is that $\beta$ is not point-wise
fixed by the reflection, in which case, there are exactly two fixed points of the reflection on it.
These two points divide the loop into two arcs which are permuted by the reflection.

\smallskip
Note that if $\beta$ is stabilized by an elliptic transformation $a \neq I$ (whose fixed points are separated by $\beta$) and an imaginary reflection $\tau$, then the transformation $b=a\tau$ must be either a reflection or an imaginary reflection (by Proposition \ref{reversing}). If $b$ is an imaginary reflection, then $\tau$ and $b$ must coincide on $\beta$, that is, they are the same, a contradiction to the fact that $a \neq I$.
Summarizing all the above is the following.

\begin{prop} \label{loopsestab}
If $\beta \in {\mathcal G}$, then $K_{\beta}$ is either:
\begin{itemize}
\item[({i})] trivial;

\item[({ii})] a cyclic group generated by an elliptic element of finite order whose fixed points are separated by $\beta$;

\item[({iii})] a cyclic group generated by an elliptic element of order two with its fixed points on $\beta$;

\item[({iv})] a cyclic group generated by an imaginary reflection;

\item[({v})] a cyclic group generated by a reflection with $\beta$ as its circle of fixed points;

\item[({vi})] a cyclic group generated by a reflection such that $\beta$ contains exactly two fixed points of it;

\item[({vii})] a group generated by a reflection, with exactly two fixed points on $\beta$, and an elliptic involution with these two points as its fixed points.
In this case, the composition of these two is a reflection with $\beta$ as circle of fixed points and  $K_{\beta}$ is isomorphic to ${\mathbb Z}_{2}^{2}$;
\item[({viii})] a group generated by an elliptic involution, whose fixed points are separated by $\beta$, and a reflection with $\beta$ as
its circle of fixed points. In this case,  the composition of these two is an imaginary reflection and $K_{\beta}$ is isomorphic to ${\mathbb Z}_{2}^{2}$;
\item[(ix)] a group generated by an elliptic involution, with both fixed points on $\beta$, and a reflection whose circle of fixed points intersects $\beta$ at two points and separating the fixed points of the elliptic involution.
In this case, $K_{\beta}$ is isomorphic to ${\mathbb Z}_{2}^{2}$; and

\item[(x)] a group generated by a reflection, with exactly two fixed points on
$\beta$, and an elliptic involution with both fixed points on the circle of fixed points of the reflection.
In this case, $K_{\beta}$ is isomorphic to ${\mathbb Z}_{2}^{2}$;
 \end{itemize}
\end{prop}

\subsection{Neighboring structural regions}
Later, we will construct a connected compact domain by gluing a finite set of structure regions an loops such that no two of these structure regions are $K$-equivalent. Let us mention the following simple fact concerning two structural regions sharing a common structure loop.

\begin{prop}\label{equivalente}
Let $R$ and $R'$ be any two different structure regions with a common boundary loop $\beta$. Then, they are
$K$-equivalent if and only if either:
\begin{itemize}
\item[({i})] $K_{\beta}$ contains an element (necessarily of order two) which does not belong
to $K_{R}$ or
\item[({ii})] there is another boundary loop $\beta'$ of $R$ and an element $k \in
K \setminus K_{R}$ so that $k(\beta')=\beta$ (necessarily a loxodromic or pseudo-hyperbolic).
\end{itemize}
\end{prop}
\begin{proof}
If two structure regions share a boundary structure loop, then they are $K$-equivalent if there is some element of $K$ sending a boundary loop of one to a boundary loop of the other. \end{proof}

\subsection{Structure regions with trivial or cyclic stabilizer}
Next result is  related to those structure regions with $K$-stabilizer being either trivial or a cyclic group generated by
a reflection.

\begin{prop}\label{trivial}
Let $R$ be a structure region with $K_{R}$ either trivial or a cyclic group generated by a reflection.
If $\beta$ a boundary loop of $R$ such that $K_{R} \cap K_{\beta}=\{I\}$, then there is
a non-trivial element $k \in K \setminus K_{R}$ so that $k(\beta)$ still a boundary loop of $R$.
\end{prop}
\begin{proof}
The hypothesis on $K_{R}$ asserts that the only possibilities in
Proposition \ref{loopsestab} for $K_{\beta}$ are ({i}), ({iii}), ({iv}) and ({v}).
In all of these cases, with the exception of ({i}), $K_{\beta}$ contains an element outside $K_{R}$.
In case ({i}) $K_{\beta}$ is trivial. The projection of $\beta$ on $S^{+}$ is a simple loop
$\beta_{*}$ which has trivial $J$-stabilizer. We have that $\beta$ is free homotopic to the product of the
other boundary loops of $R$. If none of the other boundary loops of $R$ is equivalent to $\beta$ under
$K$, then we may delete $\beta_{*}$ and its $J$-translates from ${\mathcal F}$,  contradicting the
minimality of ${\mathcal F}$.
\end{proof}

\subsection{Structure loops with non-trivial conformal stabilizer}
Let $R$ be a structure region with $K_{R}$ neither trivial or a cyclic group generated by a reflection.
Proposition \ref{regstab} asserts that $K_{R}$ is either a cyclic group generated by an elliptic transformation $a$ or it is a dihedral group generated by two reflections, say $\tau_{1}$ and $\tau_{2}$; we set $a=\tau_{2}\tau_{1}$.
In either situation, $H=K_{R}^{+}=\langle a \rangle$ (a non-trivial finite cyclic group).
The structure region  $R$ has either $0$, $1$ or $2$ boundary structure loops being stabilized by $a$. The other boundary structure loops, if any, have trivial $H$-stabilizers.

\smallskip
If no structure loop on the boundary of $R$ is stabilized by $a$, then both fixed points of $a$ lie in $R$. Also, a boundary structure loop is stabilized by $a$ if and
only if it separates the fixed points of $a$.

\smallskip
Assume now that $K_{R}$ is a dihedral group (which is known to be generated by the two reflection $\tau_{1}$ and $\tau_{2}$) and that $R$ contains only one of the fixed points of $a$. Then we have a structure loop $\beta$ in the boundary of $R$ which is invariant under $a$, so it is also invariant under $\tau_{1}$ and $\tau_{2}$. By proposition \ref{loopsestab}, $K_{R} \cong {\mathbb Z}_{2}^{2}$.

\begin{prop}\label{prop413}
 Let $R$ be a structure region with $K_{R}^{+}=H$ being non-trivial.  If there is a fixed point of $H$
in $R$, then both fixed points of $H$ lie in $R$.
\end{prop}
\begin{proof}
Suppose there is only one fixed point in $R$ of the cyclic group $H$.
Then there is a unique structure loop $\beta$ on the boundary of $R$ stabilized by $H$. Every other structure
loop, on the boundary of $R$, has $H$-stabilizer the identity.
If $K_{R}=H$, then it follows that if were to fill in the discs bounded by the other structure loops on
the boundary of $R$, then $\beta$ would be contractible; that is, if we delete the projection of $\beta$ and their
$J$-translates from our list of base structure loops, this would leave unchanged the smallest normal subgroup
containing the base structure loops raised to appropriate powers. Since we have chosen our collection ${\mathcal F}$
to be minimal, this cannot be.
If $K_{R} \neq H$, then we have a reflection $\tau \in K_{R}$ whose circle of fixed points is
not completely contained in $R$, and $K_{R}=\langle H,\tau\rangle$ is a dihedral group (as $H$ is the orientation-preserving half of $K_{R}$, it is a normal subgroup). In this case,
the fixed point of $H$ contained in $R$ is also fixed by $\tau$ (so both fixed points of $H$ are fixed
by $\tau$ as $\tau H \tau=H$), $K_{R} \cong {\mathbb Z}_{2}^{2}$ and  $H \cong {\mathbb Z}_{2}$. In this case, we
may also delete the projection of $\beta$ and its $J$-translates from ${\mathcal F}$ in order to get a contradiction
to the minimality of ${\mathcal F}$.
\end{proof}

The previous result asserts that a non-trivial elliptic transformation in $K_{R}^{+}$ either has both fixed points on the structure region $R$ or none of them belong to it. In the last case, there are (exactly) two boundary structure loops of $R$, each one invariant under such an elliptic transformation.

\begin{prop}\label{prop414}
Let $R$ be a structure region with non-trivial $K_{R}^{+}=H$.
If $\beta_{1}, \beta_{2}$ are two different boundary loops of $R$ which are invariant
under $H$, then there is a (non-trivial) element $k \in K$ so that $k(\beta_{1})=\beta_{2}$ (such an element is either loxodromic or pseudo-hyperbolic).
\end{prop}
\begin{proof}
Let us assume that there is no such element of $K$ as desired and let $R_{*}$ be the other structure region sharing $\beta_{2}$ in its boundary. Proposition \ref{prop413} asserts that on the region $R_{*}$ there is another boundary loop $\beta_{3}$ which is also invariant under $H$. All other boundary loops of $R \cup R_{*}$ (with the exception of $\beta_{1}$, $\beta_{2}$ and $\beta_{3}$) have trivial $K$-stabilizers. In particular, they are not $K$-equivalents to $\beta_{1}$, $\beta_{2}$ and $\beta_{3}$. Also, $\beta_{2}$ is neither $K$-equivalent to $\beta_{1}$ and $\beta_{3}$ (by our assumption). If we project the region $R \cup R_{*} \cup \beta_{2}$ on $S^{+}$, we obtain an homeomorphic copy and the projected loop from $\beta_{2}$ is not $J$-equivalent to none of its boundary loops. In particular, we may delete it (and its $J$-translates) obtaining a contradiction to the minimality of ${\mathcal F}$.
\end{proof}

\subsection{Structure regions with trivial stabilizers}
Let $R$ be a structure region with trivial stabilizer
$K_{R}=\{I\}$. As consequence of Proposition \ref{trivial}, every other structure region is necessarily $K$-equivalent
to $R$. It follows that $R$ is a fundamental domain for $K$ and the boundary loops are paired by
either reflections, imaginary reflections, loxodromic transformations or pseudo-hyperbolics. In this case we obtain that $K$ is the free product, by the Klein-Maskit combination theorem, of groups of types ({i})-({iv}) as described in the theorem.

\subsection{Structure regions with non-trivial stabilizers}
Let us now assume there is no structure region with trivial $K$-stabilizer.
Proposition \ref{equivalente}, and the fact that $S$ is compact and connected, permits us to construct a
connected set $\widehat{R}$ obtained as the union of a finite collection of $K$-non-equivalent structure regions
(each of them has non-trivial $K$-stabilizer) together their boundary structure loops.

\smallskip
\noindent
{\it Modification of $\widehat{R}$.}
We proceed to modify $\widehat{R}$ as follows. Let $\beta$ be a structure loop in the boundary of $\widehat{R}$ and $R \subset \widehat{R}$ be the structure region with $\beta$ in its border.
Assume there is a reflection $\tau \in K_{R}$ keeping invariant $\beta$ (so its circle of fixed points intersects $\beta$ at two points). We know that
there is some $k \in K$ such that $k(\beta)=\beta'$ is a boundary loop of $\widehat{R}$. If $\beta' \neq \beta$ (so $k$ is either loxodromic or pseudo-hyperbolic), let
$R' \subset \widehat{R}$ be the structure region containing $\beta'$ in its border. The reflection $\tau'=k \tau k^{-1}$ belongs to $K_{R'}$ and keeps invariant $\beta'$. In the case that $\tau' \neq \tau$, we eliminate $R'$ from $\widehat{R}$ and we add to it the structure region $k^{-1}(R')$. Under this type of process, we may now assume that $\tau'=\tau$.

\smallskip
The above permits us to assume that our set $\widehat{R}$ has the following extra property: {\it if $\beta$ and $\beta'$ are boundary loops of $\widehat{R}$ such that ({i}) there is some $k \in K$ with $k(\beta)=\beta'$ and ({ii}) there is a reflection $\tau \in K_{R}$ keeping invariant $\beta$, then $\tau$ also keeps invariant $\beta'$.}

\smallskip
\noindent
{\it Step 1: Internal structural loops and amalgamated free products.}
Note that if $\beta$ is a structure loop contained in the interior of $\widehat{R}$, then there are two different structure regions $R_{1}, R_{2} \subset \widehat{R}$ with $\beta$ as their common boundary. In this case, as $R_{1}$ and $R_{2}$ are non-$K$-equivalent, $K_{\beta}=K_{R_{1}} \cap K_{R_{2}}$. By Propositions \ref{prop413} and \ref{prop414} $K_{\beta}^{+}$ is trivial. So, either $K_{\beta}$ is trivial or generated by a reflection with exactly two fixed points on $\beta$.
We now perform the amalgamated free product of $K_{R_{1}}$ and $K_{R_{2}}$ with amalgamation at $K_{R_{1}} \cap K_{R_{2}}$. In the case $K_{\beta}$ is trivial, we are obtaining free product of groups of the types as described in the theorem and in the other case, we are constructing parts of the basic extended dihedral groups.

\smallskip
\noindent
{\it Step 2: Boundary structural loops and HNN-extensions.}
Next, let $\beta$ be a structure loops in the boundary of $\widehat{R}$ and let $R\subset \widehat{R}$ be the structure region with $\beta$ in its boundary.

If $K_{\beta} \cap K_{R}$ is trivial, then Proposition \ref{trivial} asserts the existence of another structure boundary loop $\beta'$ (in the boundary of $\widehat{R})$ and some $k \in K$
such that $k(\beta)=\beta'$. If $\beta' \neq \beta$, then $k$ is either loxodromic or a pseudo-hyperbolic and if $\beta'=\beta$, then $k$ is either a reflection, an imaginary reflection or an elliptic involution. In this case we obtain (by HNN extensions in the sense of the Klein-Maskit combination theorem) groups of types ({i})-({iv}) as described in the theorem.
Let us now assume $K_{\beta} \cap K_{R}$ is non-trivial.

\smallskip
\noindent
(A) Assume $K_{R}=\langle \tau \rangle$, where $\tau$ is a reflection with  circle of fixed points $C_{\tau}$. As $K_{\beta} \cap K_{R}$ is non-trivial,
$C_{\tau}$ intersects $\beta$ (at exactly two points), so either:
\begin{itemize}
\item[({i})] there is an involution $k \in K$ (conformal or anticonformal) with $k(\beta)=\beta$ and
$k(\widehat{R}) \cap \widehat{R}=\beta$; or
\item[({ii})] there is another boundary loop
$\beta'$ of $\widehat{R}$ and an element $\kappa \in K$ (which is either loxodromic or a pseudo-hyperbolic) so that $\kappa(\beta)=\beta'$ and
$\kappa(\widehat{R}) \cap \widehat{R}=\beta'$. We may assume $\kappa$ to be loxodromic. In fact, if $\kappa$ is a pseudo-hyperbolic, then $\tau\kappa$ is loxodromic with the same property.
\end{itemize}

In case ({i}), we may perform the HNN-extension (in the sense of Klein-Maskit combination theorem) to produce factors of type (2), (4) or (5) in the description of basic extended dihedral groups as in Section \ref{ejemploBEDG}.
In case ({ii}), $\beta'$ is also invariant under $\tau$. We proceed as above to produce factors of type (3) in the description of basic extended dihedral groups as in Section \ref{ejemploBEDG}.

\smallskip
The above process produces basic extended dihedral groups using the circle $C_{\tau}$ and its reflection $\tau$ (so far we have no used elliptic and loxodromics elements as in (2) or (3) in the description in Section \ref{ejemploBEDG}; these components will appear in (C) below).

\smallskip
\noindent
(B) Assume $K_{R}=\langle \phi \rangle$, where $\phi$ is an elliptic transformation. In this case we have two possibilities: either ({i}) both fixed
points of $\phi$ belong to $R$ or ({ii}) there are two boundary loops $\beta_{1}$ and $\beta_{2}$ of $R$, each one invariant under $\phi$, and there is some $k \in K$ with $k(\beta_{1})=\beta_{2}$, $k(\widehat{R}) \cap \widehat{R}=\beta_{2}$. As $K_{\beta} \cap K_{R}$ is non-trivial, necessarily
$\beta=\beta_{1}$ (or $\beta_{2}$), then (again performing HNN-extension by the Klein-Maskit combination theorem) we obtain a group of type ({vi}) as described in the theorem.

\smallskip
\noindent
(C) Assume $K_{R}$ is a dihedral group. In this case, $K^{+}_{R}$ is a non-trivial cyclic group generated by an elliptic transformation. We may proceed similarly as above to  produce basic extended dihedral groups (in this part we obtain cases (2) and (3) as described in Section \ref{ejemploBEDG}).

\smallskip
All the above asserts that $K$ is constructed as indicated  by part (1) of the theorem.
\hfill $\square$

\section{Proof of Theorem \ref{sumak=2}}\label{Sec:pruebasumak=2}
Let us consider extended Schottky groups of rank $g \geq 2$, say ${{K}}_{1}, {{K}}_{2}, {{K}}_{3}$, all of them with the same orientation-preserving half given the Schottky group $\Gamma$ of rank $g$, and set $m_{j}=m({{K}}_{j})$. All of these groups have the same region of discontinuity $\Omega$. Let us recall that each ${{K}}_{j}$ induces a symmetry $\hat{\tau}_{j}$ on the manifold $M_{\Gamma}=({\mathbb H}^{3} \cup \Omega)/\Gamma$ with exactly $m_{j}$ connected components of its set of fixed points.

\subsection{Proof of part (1)}
Let $K=\langle {{K}}_{1}, {{K}}_{2} \rangle$ and $H=K/\Gamma =\langle \hat{\tau}_{1}, \hat{\tau}_{2} \rangle \cong {\rm D}_{q}$.
Let $\theta:K \to H$ be the canonical surjection.
By Theorem \ref{const0}, the dihedral extended Schottky group $K$ is constructed using reflections $\zeta_{1}$,\ldots, $\zeta_{\alpha}$, imaginary
reflections $\eta_{1}$,\ldots, $\eta_{\beta}$, $\gamma$ cyclic groups generated by loxodromic transformations,
$\delta$ cyclic groups generated by pseudo-hyperbolics, $\rho$ cyclic groups generated by an elliptic transformation of finite order,
$\eta$ groups generated by an elliptic transformation of finite order and a loxodromic transformation commuting with it,
$\kappa$ groups generated by an elliptic transformation of finite order and a pseudo-hyperbolic (conjugating the elliptic into its inverse)
and $\varepsilon$ basic extended dihedral groups $K_{1}$,\ldots,
$K_{\varepsilon}$.

\smallskip
Each $K_{i}$ is generated by a reflection $\sigma_{i}$ and some other loxodromic and elliptic
transformations, where some of them do not commute with loxodromic ones, say $t_{i1}$,\ldots, $t_{im_{i}}$ (all of them of order at least $3$ and not commuting with $\sigma_{i}$)
and some $f_{i}$ imaginary reflections and reflections (each of them commuting with $\sigma_{i}$); let us denote
these involutions by $\sigma_{i1}$,\ldots, $\sigma_{if_{i}}$.
By the Riemann-Hurwitz formula (see Remark \ref{observa1}),
$g \geq q(\widetilde{g}-1)+1+\frac{q}{2}\sum_{i=1}^{\varepsilon}f_{i}$,
where
$\widetilde{g}=\alpha +\beta +2(\gamma+\delta+\kappa+\eta)+ \varepsilon +g_{1}+\cdots+g_{\varepsilon}-1 \geq \alpha + \beta + \varepsilon -1.
$
So, it follows that
$$
\frac{g-1}{q}+2 \geq \widetilde{g}+1+\frac{1}{2}(f_{1}+\cdots+f_{\varepsilon})\geq \alpha+\beta+\varepsilon+\frac{1}{2}(f_{1}+\cdots+f_{\varepsilon}).$$
Note that, for $\varepsilon=0$, the above also asserts that
$\left[\frac{g-1}{q}\right] \geq \alpha+\beta$.

\smallskip
A complete set of symmetries for $K$ is given by $\zeta_{1}$,\ldots, $\zeta_{\alpha}$,
$\eta_{1}$,\ldots, $\eta_{\beta}$, $\sigma_{i}$, $\sigma_{ik}$, where  $k=1,\ldots,f_{i}$  and  $i=1,\ldots,\varepsilon$.
If $c$ denotes any of the above symmetries, then $\langle c \rangle < {\rm C}
(\widehat{K},c)$ and so
$$\langle \theta(c) \rangle \subseteq \theta({\rm C} (K,c))\subseteq
{\rm C} (H,\theta(c)) \cong \left\{\begin{array}{ll}
{\rm Z}_{2}^{2}, & \mbox{$q$ even}\\
{\rm Z}_{2}, & \mbox{$q$ odd}
\end{array}
\right.
$$
Now, it is not hard to see (by looking at the structural picture of $K$) that:
$$\theta({\rm C}(K;\zeta_{j}))=\langle \theta(\zeta_{j})\rangle={\mathbb Z}_{2}, \;
\theta({\rm C}(K;\eta_{j}))=\langle \theta(\eta_{j})\rangle={\mathbb Z}_{2},\; \theta({\rm C}(K;\sigma_{jk}))=\theta(\langle \sigma_{j}, \sigma_{jk}\rangle)={\mathbb Z}_{2}^{2},$$
$$\theta({\rm C}(K;\sigma_{j}))=\theta(\langle \sigma_{j},
\sigma_{j1},\ldots,\sigma_{jf_{j}},\epsilon_{j1},\ldots, \epsilon_{jm_{j}} \rangle)={\mathbb Z}_{2}^{2}.
$$
Finally,  it follows from Theorem \ref{teo1} that
$m_{1}+m_{2} \leq 2(\alpha + \beta )+\varepsilon +(f_{1}+\cdots +f_{\varepsilon }).$
If $\varepsilon \geq 1$, then
$m_{1}+m_{2} <
2(\alpha + \beta +\varepsilon )+(f_{1}+\cdots +f_{\varepsilon}) \leq
2\left(\frac{g-1}{q}\right)+4.$
If $\varepsilon=0$, then $m_{1}+m_{2} \leq  2(\alpha+\beta) \leq
2\left[\frac{g-1}{q}\right]+4$, so we are done.
\hfill$\Box$

\subsection{Proof of part (2)}
Now, we set $\widehat{K}=\langle {{K}}_{1}, {{K}}_{2}, {{K}}_{3}\rangle$, $H=\widehat{K}/\Gamma=\langle \hat{\tau}_{1}, \hat{\tau}_{2}, \hat{\tau}_{3}\rangle$, where $\hat{\tau}_{j}$ is the involution induced from ${{K}}_{j}$, and consider the canonical surjective homomorphism $\theta:\widehat{K} \to H$ with $\Gamma$ as its kernel.

\smallskip
For each $j \in \{1,2,3\}$ and $r<s \in \{1,2,3\}-\{j\}$, we consider the dihedral extended Schottky group
$\widehat{K}_{j}=\langle {{K}}_{r},{{K}}_{s}\rangle=\theta^{-1}(\langle\hat{\tau}_{r},\hat{\tau}_{s}\rangle)$, which also have region of discontinuity $\Omega$.
Let $p_{j}=q_{r,s}$ be the order of $\hat{\tau}_{r}\hat{\tau}_{s}$; so $\langle \hat{\tau}_{r},\hat{\tau}_{s}\rangle \cong {\rm D}_{p_{j}}$ .

\smallskip
By a permutation of the indices, we may
assume that $2 \leq q_{12} \leq q_{13} \leq q_{23}$. As consequence of Part (1)  one has the
inequality
\begin{equation}\label{ec1}
m_{1}+m_{2}+m_{3} \leq \left[ \frac{g-1}{q_{12}}\right] +\left[ \frac{g-1}{q_{13}}\right]+\left[
\frac{g-1}{q_{23}}\right]+6.
\end{equation}

\subsubsection{Case $ g=2$}
In this case, $m_{1}+m_{2}+m_{3} \leq 6$. As $m_{i}+m_{j} \leq 4$, either $m_{1}+m_{2}+m_{3} \leq 5$ or else $m_{1}=m_{2}=m_{3}=2$. Claim \ref{hecho1} ends with the proof for $g=2$.
\hfill$\Box$

\begin{claim}\label{hecho1}
The case $m_{1}=m_{2}=m_{3}=2$ is not possible for $g=2$.
\end{claim}
\begin{proof}
Let $(\alpha,\beta,\gamma,\delta,\rho,\eta,\kappa,\varepsilon)$ be the signature of the dihedral extended Schottky group $\widehat{K}_{j}$.
Let $\Gamma_{1},, \ldots , \Gamma_{\varepsilon}$ be the corresponding $\varepsilon$ basic extended dihedral groups.
If $\Omega_{i}$ is the  region of discontinuity of $\Gamma_{i}$, then $\Omega_{i}/\Gamma_{i}^{+}$ has genus $g_{i}$, has
$2r_{i}$ conical points with cone orders $t_{i1},t_{i1},\ldots, t_{ir_{i}},t_{ir_{i}}$ (each of them bigger tan two) and $2n_{i}$ conical points of order two.
Then, $\Omega/\widehat{K}_{j}^{+}$ has genus
$\widetilde{g}=\alpha+\beta+2(\gamma+\delta+\kappa+\eta)+\varepsilon +g_{1}+\cdots+g_{\varepsilon} -1$,
with $2(r_{1}+\cdots+r_{\varepsilon})$ conical points of orders $t_{11}$, $t_{11}$, $t_{12}$,
$t_{12}$\ldots, $t_{\varepsilon r_{\varepsilon}}$, $t_{\varepsilon r_{\varepsilon}}$,
$2(n_{1}+\cdots+n_{\varepsilon})$ conical points of order $2$, and $4\rho$ conical points of orders $l_{1},l_{1},l_{1},l_{1},\ldots, l_{\rho},l_{\rho},l_{\rho},l_{\rho}$,
and
$$2=p_{j}(\widetilde{g}-1)+1+2p_{j}\sum_{j=1}^{\rho} \left(1-\frac{1}{l_{j}}\right)+p_{j}\sum_{i=1}^{\varepsilon}\sum_{k=1}^{r_{i}}\left(1-\frac{1}{t_{ik}}\right)+
\frac{p_{j}}{2}\sum_{i=1}^{\varepsilon}n_{i} \geq p_{j}(\widetilde{g}-1)+1.
$$
It follows that $\widetilde{g} \in \{0,1\}$.

\smallskip
If $\alpha=\beta=\varepsilon=0$ (so $\delta>0$), then $\widetilde{g}=2(\gamma+\delta+\kappa+\eta)-1$, so $\widetilde{g}=1=\delta+\kappa$ and $\eta=\gamma=0$. It follows that $\widehat{K}_{j}$ is a elementary group (either ({i}) generated by a pseudo-hyperbolic or ({ii}) a pseudo-hyperbolic and an elliptic element with the same fixed points), a contradiction (as it must contain a Schottky group of rank two).

\smallskip
If $\alpha+\beta+\varepsilon>0$, it must be that $\gamma=\delta=\kappa=\eta=0$ and $\widetilde{g}=\alpha+\beta+\varepsilon+g_{1}+\cdots+g_{\varepsilon}-1$. If $\varepsilon=0$, then $\alpha+\beta \in \{1,2\}$ and $\widehat{K}_{j}$ will be an
elementary group, again a contradiction. So $\varepsilon \geq 1$. As
$\alpha+\beta+\varepsilon +g_{1}+\cdots+g_{\varepsilon} =\widetilde{g}+1\in \{1,2\}$,
we must have $\varepsilon=1$. In fact, if $\varepsilon \geq 2$, then $\varepsilon=2$, $\alpha=\beta=g_{1}=g_{2}=0$, $\widetilde{g}=1$ and
$$1=p_{j}\sum_{i=1}^{2}\sum_{k=1}^{r_{i}}\left(1-\frac{1}{t_{ik}}\right)+
\frac{p_{j}}{2}\sum_{i=1}^{2}n_{i}.$$
If, for instance, $n_{1} \geq 1$, then (as $p_{j} \geq 2$) $r_{1}=r_{2}=0=n_{2}$ and $n_{1}=1$, from which we obtain that $\widehat{K}_{j}$ is elementary, a contradiction. In particular, $n_{1}=n_{2}=0$. The other case can be worked similarly to obtain a contradiction.

\smallskip
Assume now that $\varepsilon=1$. Then $\alpha+\beta+g_{1} \in \{0,1\}$. If $g_{1}=1$, then $\alpha=\beta=0$ and $\widehat{K}_{j}$ is again elementary. So, $g_{1}=0$
and $\widetilde{g}=\alpha+\beta \in \{0,1\}$, that is
$(\alpha,\beta) \in \{(0,0), (1,0), (0,1)\}.$
In the case $(\alpha,\beta)=(0,0)$, we have $\widetilde{g}=0$ and $$2 \geq 1-n
+n\sum_{k=1}^{r_{1}}\left(1-\frac{1}{t_{ik}}\right)+ \frac{nn_{1}}{2}.$$
As each $|t_{ij}| \geq 2$, the above ensures that
$2 \geq n(n_{1}+r_{1}-2).$
Now, if $n_{1}+r_{1} \geq 3$, then $n=2$ (as required) and, necessarily, $n_{1}+r_{1}=3$. If $n_{1}+r_{1} \leq
2$, then $\widehat{K}_{j}$ will be elementary, a contradiction.
In the case $(\alpha,\beta) \in \{(1,0),(0,1)\}$, we have $\widetilde{g}=1$
and $$2 \geq 1+n\sum_{k=1}^{r_{1}}\left(1-\frac{1}{t_{ik}}\right)+ \frac{nn_{1}}{2}.
$$
Again, as $|t_{1k}| \geq 2$, one obtains that
$2 \geq n(n_{1}+r_{1}).$ Then $n=2$ and $n_{1}+r_{1}=1$, but again this makes $\widehat{K}_{j}$ to be elementary, a contradiction.

\smallskip
As a consequence of the above, we see that each $\widehat{K}_{j}$ is generated by a reflection $\rho_{j1}$ and three other
involutions $\rho_{j2}$, $\rho_{j3}$, $\rho_{j4}$, each one commuting with $\rho_{j1}$ (some of the involutions
may be reflections and others may be imaginary reflections). Note also that $\widehat{K}_{j}$ keeps invariant the circle of
fixed points of $\rho_{j1}$, that is, the limit set of $\widehat{K}_{j}={\mathbb Z}_{2} \oplus ({\mathbb Z}_{2}*{\mathbb
Z}_{2}*{\mathbb Z}_{2})$ is contained in the circle of fixed points of $\rho_{j1}$. As $\widehat{K}_{j}$ has finite index in $\widehat{K}$, the limit set of $\widehat{K}$ (which is the same as for $\widehat{K}_{j}$) is contained in such a circle. Since
the limit set of
$\widehat{K}$ is infinite, such a circle is uniquely determined by it. As a consequence,
$\rho_{11}=\rho_{21}=\rho_{31}=\rho$, in particular, the symmetry $\theta(\rho) \in \theta(\widehat{K}_{1}) \cap \theta(\widehat{K}_{2}) \cap
\theta(\widehat{K}_{3})$. As the symmetries in $H$ are exactly $\hat{\tau}_{1}$, $\hat{\tau}_{2}$,
$\hat{\tau}_{3}$ and $\hat{\tau}_{1}\hat{\tau}_{2}\hat{\tau}_{3}$ and, by the definition of the $\widehat{K}_{j}$'s,
$\hat{\tau}_{j} \notin \theta(\widehat{K}_{1}) \cap \theta(\widehat{K}_{2})
\cap \theta(\widehat{K}_{3})$, we must have
$\theta(\rho)=\hat{\tau}_{1}\hat{\tau}_{2}\hat{\tau}_{3}$. In particular, $\hat{\tau}_{1}\hat{\tau}_{2}\hat{\tau}_{3} \in \langle
\hat{\tau}_{1},\hat{\tau}_{2}\rangle$, obligating to have that $\hat{\tau}_{3} \in \langle \hat{\tau}_{1},\hat{\tau}_{2}\rangle$, a
contradiction.
\end{proof}

\subsubsection{Case $g=3$}
If $q_{23} \geq 3$, then $m_{2}+m_{3} \leq 4$. Since $m_{1} \leq g+1=4$, we get in this case that
$m_{1}+m_{2}+m_{3} \leq 8$.
Let us now assume $q_{12}=q_{13}=q_{23}=2$, in which case $H=\langle \hat{\tau}_{1},\hat{\tau}_{2},\hat{\tau}_{3}\rangle \cong
{\mathbb Z}_{2}^{3}$. We may reorder again the indices to assume that $m_{1} \leq m_{2} \leq m_{3} \leq 4$. If
$m_{3}=4$, then inequality \eqref{ec1} asserts that $m_{1}+m_{2} \leq 2$, so $m_{1}+m_{2}+m_{3} \leq 6$.  Now,
the only case with $m_{3} \leq 3$ for which $m_{1}+m_{2}+m_{3} > 8$ is when
$m_{1}=m_{2}=m_{3}=3$. Claim \ref{claim5.2} ends the proof for $g=3$.
\hfill$\Box$

\begin{claim}\label{claim5.2}
The case $m_{1}=m_{2}=m_{3}=3$ is not possible for $g=3$.
\end{claim}
\begin{proof}The Proof follows the same ideas as for the proof of Claim \ref{hecho1}. If we assume $m_{1}=m_{2}=m_{3}=3$, then
 $H=\langle \hat{\tau}_{1},\hat{\tau}_{2},\hat{\tau}_{3}\rangle={\mathbb Z}_{2}^{3}$.

\smallskip
As before, we consider the dihedral extended Schottky groups $\widehat{K}_{j}=\theta^{-1}(\langle
\hat{\tau}_{r},\hat{\tau}_{s}\rangle)$, where $r<s$ and $r,s \in \{1,2,3\}-\{j\}$.
So, $\theta(\widehat{K}_{j})={\mathbb Z}_{2}^{2}$ and $\widehat{K}_{j}$ has index two in $\widehat{K}$.

\smallskip
Let $(\alpha,\beta,\gamma,\delta,\rho,\eta,\kappa,\varepsilon)$ be the signature of the dihedral extended Schottky group
$\widehat{K}_{j}$, and let $\Gamma_{1},  \ldots, \Gamma_{\varepsilon}$ be the correspondent $\varepsilon$ basic extended dihedral groups. Recall that
$\alpha+\beta+\delta +\varepsilon>0$.  As before, $\Omega/\widehat{K}_{j}^{+}$ will have genus
$\widetilde{g}=\alpha+\beta+2(\gamma+\delta+\kappa+\eta)+\varepsilon+g_{1}+\cdots+g_{\varepsilon} -1,$
with $2(2\rho+r_{1}+\cdots+r_{\varepsilon}+n_{1}+\cdots+n_{\varepsilon})$ conical points of order $2$,  and so that
$$2 \geq 2(\widetilde{g}-1)+\sum_{i=1}^{\varepsilon}(r_{i}+n_{i}).$$
In particular, $\widetilde{g} \in \{0,1\}$.
If $\alpha=\beta=\varepsilon=0$ (so $\delta+\kappa>0$), then $\widetilde{g}=2(\gamma+\delta+\kappa+\eta)-1$, so $\gamma=\eta=0$ and $\delta+\kappa=1$, from which we obtain
that $\widehat{K}_{j}$ is elementary, a contradiction.
As $\alpha+\beta+\varepsilon>0$, then $\gamma=\delta=\kappa=\eta=0$. If $\varepsilon=0$, the same
asserts that $\widehat{K}_{j}$ is an elementary group containing the non-elementary group $\Gamma$, a contradiction, so
$\varepsilon \geq 1$. If $\varepsilon \geq 3$, then $\alpha+\beta<0$, a contradiction.
If $\varepsilon=2$, then $\alpha+\beta=\widetilde{g}-1$, so the only possible case is to
have $\widetilde{g}=1$ and $\alpha=\beta=0$. In this way,
$$(\alpha,\beta,\varepsilon;\widetilde{g}) \in
\{(0,0,1;0), (1,0,1;1), (0,1,1;1), (0,0,2;1)\}.$$
Let us consider the case $(\alpha,\beta,\varepsilon;\widetilde{g})=(0,0,2;1)$. In this case
$\widehat{K}_{j}$ is free product of two groups $\Gamma_{1}=\langle \rho_{1},\eta_{1}\rangle={\mathbb Z}_{2}^{2}$ and
$\Gamma_{2}=\langle \rho_{2},\eta_{2}\rangle={\mathbb Z}_{2}^{2}$, where $\rho_{k}$ is reflection and $\eta_{k}$
is either a reflection or an imaginary reflection. In this case a complete set of symmetries of $\widehat{K}_{j}$ is given
by $\{\rho_{1},\rho_{2},\eta_{1},\eta_{2}\}$. Let us note that $\theta(\rho_{k}) \neq \theta(\eta_{k})$ since
otherwise the elliptic element of order two $\rho_{k}\eta_{k} \in \Gamma$, a contradiction to the fact that
$\Gamma$ is torsion free. As
$${\rm C}(\widehat{K}_{j},\rho_{k})={\rm C}(\widehat{K}_{j},\eta_{k})=\langle
\rho_{k},\eta_{k}\rangle={\mathbb Z}_{2}^{2}$$ and that $\theta(\rho_{k}) \neq \theta(\eta_{k})$, we may see from
Theorem \ref{teo1} that the number of components of fixed points of any of the symmetries $\theta(\rho_{1})$ and
$\theta(\eta_{1})$ (in the manifold $M_{\Gamma}$) is at most $2$, a contradiction to the fact we are assuming $m_{1}=m_{2}=m_{3}=3$.

\smallskip
Then, let us consider the case $(\alpha,\beta,\varepsilon; \widetilde{g})=(0,0,1;0)$. In this case
$\widehat{K}_{j}=\langle \rho, \eta_{1},\eta_{2},\eta_{3},\eta_{4}\rangle$, where $\rho$ is a reflection and each of the
$\eta_{k}$ is either a reflection or an imaginary reflection commuting with $\rho$ (each element of $\widehat{K}_{j}$ keeps
invariant the circle of fixed points of $\rho$). As before, $\theta(\rho) \neq \theta(\eta_{k})$, for each
$k=1,2,3,4$. So, in particular, $\theta(\eta_{1})=\theta(\eta_{2})=\theta(\eta_{3})=\theta(\eta_{4})$. In this
case, a complete set of symmetries of $\widehat{K}_{j}$ is given by $\{\rho, \eta_{1},\eta_{2},\eta_{3},\eta_{4}\}$ and
$${\rm C}(\widehat{K}_{j},\rho)=\widehat{K}_{j}, \; {\rm C}(\widehat{K}_{j},\eta_{k})=\langle \rho, \eta_{k}\rangle={\mathbb Z}_{2}^{2}$$
and it follows, from Theorem \ref{teo1}, that $\theta(\rho)$ has at most one connected components
of fixed points on $M_{\Gamma}$, a contradiction to the assumption that each of the symmetries have $3$ connected components.

\smallskip
Let now $(\alpha,\beta,\varepsilon; \widetilde{g})=(1,0,1;0)$ (similar
arguments for the case $(\alpha,\beta,\varepsilon;\widetilde{g})=(0,1,1;0)$). In this case
$\widehat{K}_{j}$ is a free product of a cyclic group generated by a reflection $\zeta$ and a group $\widehat{K}_{j}^{0}= \langle
\rho, \eta_{1},\eta_{2},\eta_{3},\eta_{4}\rangle$, where $\rho$ is a reflection and each of the $\eta_{k}$ is
either a reflection or an imaginary reflection commuting with $\rho$ (each element of $\widehat{K}_{j}^{0}$ keeps invariant
the circle of fixed points of $\rho$). Again, as before, $\theta(\rho) \neq \theta(\eta_{k})$, for each
$k=1,2,3,4$. So, in particular, $\theta(\eta_{1})=\theta(\eta_{2})=\theta(\eta_{3})=\theta(\eta_{4})$. In this
case, a complete set of symmetries of $\widehat{K}_{j}$ is given by $\{\zeta, \rho, \eta_{1},\eta_{2},\eta_{3},\eta_{4}\}$
and
$${\rm C}(\widehat{K}_{j},\zeta)=\langle \zeta\rangle={\mathbb Z}_{2}, \; {\rm C}(\widehat{K}_{j},
\rho)=\widehat{K}_{j}^{0}, \; {\rm C}(\widehat{K}_{j},\eta_{k})=\langle \rho, \eta_{k}\rangle={\mathbb Z}_{2}^{2}.$$
The only way for both symmetries in $\theta(\widehat{K}_{j})$ to have exactly $3$ connected components of fixed points on $M$ is
to have that $\theta(\zeta)=\theta(\rho)$. But in this case we should have that $\Gamma$ is the Schottky group
generated by the loxodromic transformations
$\rho\zeta, \eta_{4}\eta_{1}, \eta_{4}\eta_{2}, \eta_{4}\eta_{3},$
which is of rank $4$, a contradiction.
\end{proof}

\subsubsection{Case $g\geq 4$}
In this case, inequality \eqref{ec1} asserts
\begin{equation}\label{ec2}
m_{1}+m_{2}+m_{3} \leq \left( \frac{1}{q_{12}}+\frac{1}{q_{13}}+\frac{1}{q_{23}}\right)(g-1)+6.
\end{equation}
If $q_{12}^{-1}+q_{13}^{-1}+q_{23}^{-1}\leq 1$, then the above ensures that
$m_{1}+m_{2}+m_{3} \leq g+5$.
If $q_{12}^{-1}+q_{13}^{-1}+q_{23}^{-1}> 1$, then we have the following cases:
\begin{equation}\label{tablita}
\begin{array}{|c|c|c|c|}\hline
&&&\\[-3mm]
\hspace{4mm} q_{12} \hspace{4mm}  & \hspace{4mm}  q_{13} \hspace{4mm}  & \hspace{4mm}  q_{23}\hspace{4mm}  &  \hspace{4mm}  H=\langle \hat{\tau}_{1},\hat{\tau}_{2},\hat{\tau}_{3}\rangle \hspace{4mm}  \\
&&&\\[-3mm]
\hline
2& 2 & r \geq 2 & {\mathbb Z}_{2} \times {\rm D}_{r} \\
2 & 3 & 3 & {\mathbb Z}_{2} \ltimes {\rm A}_{4} \\
2 & 3 & 4 & {\mathbb Z}_{2} \ltimes {\rm  S}_{4} \\
2 & 3 & 5 & {\mathbb Z}_{2} \ltimes {\rm  A}_{5} \\\hline
\end{array}
\end{equation}

\medskip \noindent
where ${\mathbb Z}_{2}=\langle \hat{\tau}_{1}\rangle$, ${r}=\langle \hat{\tau}_{2},\hat{\tau}_{3}\rangle$ and ${{\rm A}}_{4}$,
${\rm  S}_{4}$ and ${{\rm A}}_{5}$ are generated by $\langle \hat{\tau}_{1}\hat{\tau}_{2}, \hat{\tau}_{1}\hat{\tau}_{3}\rangle$
in each case. In the cases $q_{13}=3$ one has that $\langle \hat{\tau}_{1},\hat{\tau}_{3}\rangle \cong {\rm D}_3    $, so $\hat{\tau}_{1}$ and
$\hat{\tau}_{3}$ are conjugated, that is, $m_{1}=m_{3}$. It follows that
$$2m_{1}=2m_{3}=m_{1}+m_{3} \leq 2 \left[ \frac{g-1}{3}\right]+4,$$
so $\displaystyle{m_{1}=m_{3} \leq \left[ \frac{g-1}{3}\right]+2}$. As $q_{23} \geq 3$,
$$m_{1}+(m_{2}+m_{3}) \leq \left(\left[ \frac{g-1}{3}\right]+2\right) +
\left( 2 \left[ \frac{g-1}{q_{23}}\right]+4\right)\leq 3\left[ \frac{g-1}{3}\right]+6 \leq g+5.$$
Now, for the case $q_{12}=q_{13}=2$ and $q_{23}=r \geq 2$, inequality \eqref{ec1} asserts
$$m_{1}+m_{2}+m_{3} \leq 2\left[ \frac{g-1}{2}\right]+\left[ \frac{g-1}{r}\right]+6 \leq 2
\left( \frac{g-1}{2}\right)+\frac{g-1}{r}+6=\frac{(r+1)g+5r-1}{r}.$$

\section{Examples}\label{Sec:4}

\begin{example}[Sharpness of the upper  bound in part (1) of Theorem \ref{sumak=2}] \rm
Fix an integer $q \geq 2$ and let ${{K}}$ be an extended Schottky group constructed (by Theorem
\ref{maintheo}) using exactly $r+1$ reflections $E_{1}$,\ldots , $E_{r+1}$. The orbifold uniformized
by  ${{K}}$ is a planar surface bounded by exactly $r+1$ boundary loops. Let us consider the surjective
homomorphism
$$\theta:{{K}} \to {\rm D}_{q}=\langle \hat{\tau}_{1},\hat{\tau}_{2}: \hat{\tau}_{1}^{2}=\hat{\tau}_{2}^{2}=(\hat{\tau}_{2}\hat{\tau}_{1})^{q}=1\rangle: \;
\theta(E_{1})=\hat{\tau}_{1}, \; \theta(E_{2})=\cdots=\theta(E_{r+1})=\hat{\tau}_{2},$$
and let $\Gamma=\ker(\theta)$ be its kernel. If we set $L=E_{2}E_{1}$ and $C_{j}=E_{r+1}E_{j+1}$, for $j=1,\ldots ,r-1$, then it is
not difficult to see that
$$\Gamma=\langle L^{q},C_{1},\ldots ,C_{r-1}, LC_{1}L^{-1},\ldots ,LC_{r-1}L^{-1},\ldots , L^{q-1}C_{1}L^{-q+1},
\ldots ,L^{q-1}C_{r-1}L^{-q+1} \rangle$$
is a Schottky group of rank $g=(r-1)q+1$. Let us consider the extended Schottky groups
$${{K}}_{1}=\theta^{-1}(\langle \hat{\tau}_{1}\rangle)=\langle E_{1}, \Gamma\rangle, \quad
{{K}}_{2}=\theta^{-1}(\langle \hat{\tau}_{2}\rangle)=\langle E_{2}, \Gamma\rangle.$$
We claim that these two extended Schottky groups  provide the upper bound in part (1) of Theorem \ref{sumak=2}.

\smallskip
As the group ${{K}}$ contains no imaginary reflections nor real Schottky groups, it follows that
${{K}}_{j}$ is constructed, by Theorem \ref{maintheo}, using $\alpha_{i}$ reflections and $\beta_{i}$ loxodromic
and pseudo-hyperbolic transformations. The handlebody $M_{\Gamma}$ admits two symmetries, say $\hat{\tau}_{1}$ and $\hat{\tau}_{2}$
induced by ${{K}}_{1}$ and ${{K}}_{2}$ respectively. The number of connected components of fixed points of $\hat{\tau}_{i}$ is
exactly $m_{i}=\alpha_{i}$. By Theorem \ref{sumak=2}, we should have $\alpha_{1}+\alpha_{2} \leq
2(r+1)$. Next, we proceed to see that in fact we have an equality, showing that the upper bound in Theorem
\ref{sumak=2} is sharp. In order to achieve the above, we use Theorem \ref{teo1}. A complete set of symmetries in
${{K}}$ is provided by $E_{1}$,\ldots , $E_{r+1}$. We also should note that ${\rm C}
(K,E_i)=\langle E_{i}\rangle$, for every $j$, and that $J(i)=\emptyset$ and $I(i)=F(i)$.

\subsubsection*{Case $q$ odd}
In this case, $I(1)=\{1,2,\ldots ,r+1\}$ and ${\rm C} ({\rm D}_{q},\hat{\tau}_{j})=\langle \hat{\tau}_{j}\rangle$. It follows, from Theorem \ref{teo1}, that
$\alpha_{1}=\alpha_{2}=r+1$
and we are done.

\subsubsection*{Case $q$ even}
In this case, $I(1)=\{1\}$, $I(2)=\{2,\ldots ,r+1\}$, ${\rm C} ({\rm D}_{q},\hat{\tau}_{1})=\langle \hat{\tau}_{1}, (\hat{\tau}_{2}\hat{\tau}_{1})^{q/2} \rangle\cong {\mathbb Z}_{2}^{2}$ and ${\rm C}
({\rm D}_{q},\hat{\tau}_{2})=\langle \hat{\tau}_{2}, (\hat{\tau}_{2}\hat{\tau}_{1})^{q/2} \rangle\cong {\mathbb Z}_{2}^{2}$. It
follows, from Theorem \ref{teo1}, that
$\alpha_{1}=2, \; \alpha_{2}=2r$
and we are done.
    \hfill $\square$ \end{example}

\begin{example}[Sharpness of the upper  bounds in part (2) of Theorem \ref{sumak=2} for $g=2$] \rm
Let ${{K}}$ be the extended Kleinian group generated by four reflections, say $\eta_{1}$,\ldots,
$\eta_{4}$, where $\eta_{1}(z)=\overline{z}$, $\eta_{2}(z)=-\overline{z}$, $\eta_{4}(z)=1/\overline{z}$ and
$\eta_{3}$ is the reflection on a circle $\Sigma$ which is orthogonal to the unit circle and disjoint from the
real and imaginary axis. One may see that
$${{K}}=\langle \eta_{4} \rangle \times ( \langle \eta_{1},\eta_{2} \rangle \ast
\langle \eta_{3} \rangle) \cong {\mathbb Z}_{2} \times ({\rm D}_{2} * {\mathbb Z}_{2}).$$
If $\Gamma=\langle A=\eta_{1}\eta_{3}, B=\eta_{1}\eta_{2}\eta_{3}\eta_{2}\rangle$, then it is a Schottky group of rank $2$
with a fundamental domain bounded by the $4$ circles
$C_{1}=\Sigma$, $C'_{1}=\eta_{1}(\Sigma)$, $C_{2}=\eta_{2}(\Sigma)$ and $C'_{2}=\eta_{2}\eta_{1}(\Sigma)$ such that
$A(C_{1})=C'_{1}$ and $B(C_{2})=C'_{2}$.
As $\eta_{1} A \eta_{1}=\eta_{3} A \eta_{3}=A^{-1}$, $\eta_{2}A\eta_{2}=\eta_{4}B\eta_{4}=B$,
$\eta_{4}A\eta_{4}=\eta_{2}B\eta_{2}=A$, $\eta_{1}B\eta_{1}=B^{-1}$ and $\eta_{3}B\eta_{3}=A^{-1}B^{-1}A$, it
follows that $\Gamma$ is a normal subgroup of ${{K}}$. Moreover, ${{K}}/\Gamma \cong {\mathbb
Z}_{2}^{3}={\mathbb Z}_{2} \times {\rm D}_{2}.$
On the handlebody $M_{\Gamma}$ both $\eta_{1}$ and $\eta_{3}$ induce the same symmetry
$\hat{\tau}_{1}$ with exactly $3$ connected components of fixed points (each of them a disc), $\eta_{2}$ induces a
symmetry $\hat{\tau}_{2}$ with exactly one connected component of fixed points (a dividing disc) and $\eta_{4}$ induces
a symmetry $\hat{\tau}_{3}$ with exactly one connected component of fixed points (this being an sphere with three
borders). It follows that the three induced symmetries are non-conjugated. So, in this case, the three extended Schottky groups ${{K}}_{j}=\theta^{-1}(\langle \hat{\tau}_{j}\rangle)$, where $\theta:{{K}} \to {{K}}/\Gamma$ is the canonical surjective homomorphism, are the ones providing the upper bound in part (2) of Theorem \ref{sumak=2} for $g=2$.
    \hfill $\square$
\end{example}

\begin{example}[Sharpness of the upper  bounds in part (2) of Theorem \ref{sumak=2} for $g = 3$] \rm
Let ${{K}}$ be the extended Kleinian group generated by three reflections,
$\eta_{1}(z)=\overline{z}$, $\eta_{2}(z)=-\overline{z}$ and $\eta_{3}$ the reflection on a circle $\Sigma$
disjoint from the real and imaginary lines. One has that
$$
{{K}}=\langle \eta_{1}, \eta_{2} \rangle * \langle \eta_{3} \rangle \cong {\rm D}_2 * {\mathbb Z}_{2}.
$$
The group $\Gamma=\langle A_{1}=(\eta_{3}\eta_{1})^{2},A_{2}=(\eta_{3}\eta_{2})^{2},
A_{3}=\eta_{3}\eta_{2}\eta_{1}\eta_{3}\eta_{1}\eta_{2}\rangle$ is a
Schottky group of rank $3$ with a fundamental domain bounded by the $6$ circles $C_{1}=\eta_{1}(\Sigma)$,
$C'_{1}=\eta_{3}(C_{1})$, $C_{2}=\eta_{2}(\Sigma)$, $C'_{2}=\eta_{3}(C_{2})$, $C_{3}=\eta_{2}(C_{1})$ and
$C'_{3}=\eta_{3}(C_{3})$,  such that $A_{1}(C_{1})=C'_{1}$, $A_{2}(C_{2})=C'_{2}$ and $A_{3}(C_{3})=C'_{3}$.
Similarly to the previous case, one may check that $\Gamma$ is a normal subgroup of ${{K}}$ and that
${{K}}/\Gamma \cong {\mathbb Z}_{2}^{3}$. The reflection $\eta_{j}$ induces a symmetry $\hat{\tau}_{j}$ (for each
$j=1,2,3$) on the handlebody $M_{\Gamma}$. In this case (either by direct inspection or by  using
Theorem \ref{teo1}) each of $\hat{\tau}_{1}$ and $\hat{\tau}_{2}$ has exactly $2$ connected components of fixed
points, and $\hat{\tau}_{3}$ has $4$ connect components of fixed points. In some cases the handlebody $M_{\Gamma}$ will have
extra automorphisms conjugating $\hat{\tau}_{1}$ with $\hat{\tau}_{2}$ (for instance, when $\Sigma$ is orthogonal to the line
$L=\{Re(z)=Im(z) \}$); but in the generic case this will not happen (that is, the three of them will be
non-conjugated).
    \hfill $\square$
\end{example}

\begin{example}[Sharpness of the upper  bounds in part (2) of Theorem \ref{sumak=2} for $g \geq 4 $ and $H \ncong {\mathbb Z}_{2} \times {\rm D}_r $] \rm
Let ${{K}}$ be an extended Schottky group constructed by using $2n+3$ reflections (or imaginary reflections
or combination of them), say $\eta_{1}$,\ldots, $\eta_{2n+3}$. Consider the surjective homomorphism
$\theta:{{K}} \to {\rm D}_{3} =\langle a,b: a^{2}=b^{2}=(ab)^{3}=1\rangle$, defined by
$\theta(\eta_{1})=\cdots =\theta(\eta_{2n+2})=a$, $\theta(\eta_{2n+3})=b$.
In this case $\Gamma=\ker(\theta)$ is a Schottky group of rank $g=6n+4$.  The handlebody $M_{\Gamma}$
admits the symmetries $\hat{\tau}_{1}=a$, $\hat{\tau}_{2}=b$ and $\hat{\tau}_{3}=bab$. By direct inspection (or
by using Theorem \ref{teo1}) it can be checked that $m_{1}=m_{2}=m_{3}=2n+3$.
    \hfill $\square$
\end{example}

\begin{example}[Sharpness of the upper  bounds in part (2) of Theorem \ref{sumak=2} for $g \geq 4$ and $H \cong {\mathbb Z}_{2} \times {\rm D}_r$] \rm
Let ${{K}}$ be an extended Schottky group constructed by using $3$ reflections (or imaginary reflections or
combination of them), say $\eta_{1}$, $\eta_{2}$ and $\eta_{3}$. Consider the surjective
homomorphism
$\theta:{{K}} \to {\mathbb Z}_{2} \times {\rm D}_3  =\langle c\rangle \times \langle a,b:
a^{2}=b^{2}=(ab)^{3}=1\rangle$, defined by
$\theta(\eta_{1})=c$, $\theta(\eta_{2})=a$ and $\theta(\eta_{3})=b$.
In this case $\Gamma=\ker(\theta)$ is a Schottky group of rank $g=2r+1$.  The handlebody $M_{\Gamma}$
admits the symmetries $\hat{\tau}_{1}=c$, $\hat{\tau}_{2}=a$ and $\hat{\tau}_{3}=b$. In this case, $m_{1}=2r$ and $m_{2}=m_{3}=4$.
    \hfill $\square$
\end{example}

\subsection*{Acknowledgments} The authors are very grateful to the reviewer for the comments and suggestions.


\end{document}